\documentclass[12pt]{article}
\input{amssym.def}
\input{amssym.tex}
\parskip=\medskipamount
\newtheorem{definition}{Definition}[section]
\newtheorem{lemma}[definition]{Lemma}
\newtheorem{theorem}[definition]{Theorem}
\newtheorem{proposition}[definition]{Proposition}

\newtheorem{remark}[definition]{Remark}

\newtheorem{examples}[definition]{Examples}

 1   %Caracteres g"ticos.
\font\ddpp=msbm10  scaled \magstep 1  %Caracteres "doble palo".

\def\QED{\hskip0.1em\hfill\null\ \null\nobreak\hfill
\kern3pt\lower1.8pt\vbox{\hrule\hbox
{\vrule\kern1pt\vbox{\kern1.7pt \hbox{$\scriptstyle
QED$}\kern0.2pt}\kern1pt\vrule}\hrule}}
\def\R{\hbox{\ddpp R}}               %Numeros reales
\def\lcf{\lbrack\! \lbrack}
\def\rcf{\rbrack\! \rbrack}

\newcommand\prueba {\mbox{{\em Proof: }}}

\begin{document}
\baselineskip=.55cm
\title{{\bf LIE ALGEBROID FOLIATIONS AND ${\cal E}^1(M)$-DIRAC
STRUCTURES}}
\author{David Iglesias, Juan C. Marrero
\\ {\small\it Departamento de Matem\'atica Fundamental, Facultad de
Matem\'aticas,}\\[-8pt] {\small\it Universidad de la Laguna, La
Laguna,} \\[-8pt] {\small\it Tenerife, Canary Islands,
SPAIN,}\\[-8pt] {\small\it E-mail: diglesia@ull.es,
jcmarrer@ull.es} }
\date{}

\maketitle \baselineskip=.45cm
\begin{abstract}
{\small We prove some general results about the relation between
the 1-cocycles of an arbitrary Lie algebroid $A$ over $M$ and the
leaves of the Lie algebroid foliation on $M$ associated with $A$.
Using these results, we show that a ${\cal E}^1(M)$-Dirac
structure $L$ induces on every leaf $F$ of its characteristic
foliation a ${\cal E}^1(F)$-Dirac structure $L_F$, which comes
from a precontact structure or from a locally conformal
presymplectic structure on $F$. In addition, we prove that a Dirac
structure $\tilde{L}$ on $M\times \R$ can be obtained from $L$ and
we discuss the relation between the leaves of the characteristic
foliations of $L$ and $\tilde{L}$. }
\end{abstract}
\begin{quote}

\end{quote}
\baselineskip=.55cm

\vspace{-1cm}
\section{Introduction}
It is well-known the fundamental role that Poisson algebras play
in Dirac's theory of constrained Hamiltonian systems \cite{Di}.
Two natural ways for Poisson algebras to arise from a manifold $M$
are through Poisson or presymplectic structures on $M$. Both
structures are examples of Dirac structures in the sense of
Courant-Weinstein \cite{Co,CW}. A Dirac structure on a manifold
$M$ is a vector sub-bundle $\tilde{L}$ of $TM\oplus T^\ast M$ that
is maximally isotropic under the natural symmetric pairing on
$TM\oplus T^\ast M$ and such that the space of sections of
$\tilde{L}$, $\Gamma (\tilde{L})$, is closed under the Courant
bracket $\makebox{{\bf [}}\, ,\, \makebox{{\bf ]}}\,^{\tilde{ }}$
on $\Gamma (TM\oplus T^\ast M)$ (see Section \ref{ejemplos},
Example 1). If $\tilde{L}$ is a Dirac structure on $M$, then
$\tilde{L}$ is endowed with a Lie algebroid structure over $M$ and
the leaves of the induced Lie algebroid foliation ${\cal
F}_{\tilde{L}}$ on $M$ are presymplectic manifolds (see
\cite{Co}). In the particular case when the Dirac structure
$\tilde{L}$ comes from a Poisson structure $\Pi$ on $M$, then
$\tilde{L}$ is isomorphic to the cotangent Lie algebroid
associated with $\Pi$ and ${\cal F}_{\tilde{L}}$ is just the
symplectic foliation of $M$ (see \cite{Co}).

An algebraic treatment of Dirac structures was developed by
Dorfman in \cite{Do} using the notion of a complex over a Lie
algebra. This treatment was applied to the study of general
Hamiltonian structures and their role in integrability. More
recently, the properties of the Courant bracket $\makebox{{\bf
[}}\, ,\, \makebox{{\bf ]}}\,^{\tilde{ }}$ have been systematized
by Liu, Weinstein and Xu \cite{LWX} in the definition of a Courant
algebroid structure on a vector bundle $E\to M$ (see also
\cite{LWX2,RW}). The natural example of a Courant algebroid is the
Whitney sum $E=A\oplus A^\ast$, where the pair $(A,A^\ast)$ is a
Lie bialgebroid over $M$ in the sense of Mackenzie-Xu \cite{MX}.

On the other hand, a Jacobi structure on a manifold $M$ is a local
Lie algebra structure, in the sense of Kirillov \cite{Ki}, on the
space $C^\infty (M,\R)$ (see \cite{DLM,GL,Li2}). Apart from
Poisson manifolds, interesting examples of Jacobi manifolds are
contact and locally conformal symplectic manifolds. In fact, a
Jacobi structure on $M$ defines a generalized foliation, the
characteristic foliation of $M$, whose leaves are contact or
locally conformal symplectic manifolds \cite{GL,Ki}. Moreover, the
1-jet bundle $T^\ast M\times \R\to M$ is a Lie algebroid and the
corresponding Lie algebroid foliation is just the characteristic
foliation of $M$ (see \cite{KS}). However, for a Jacobi manifold
$M$ the vector bundle $T^\ast M$ is not, in general, a Lie
algebroid and, in addition, if one considers the usual Lie
algebroid structure on $TM\times\R$ then the pair $(TM\times \R,
T^\ast M\times \R)$ is not a Lie bialgebroid (see \cite{IM2,V2}).
Thus, it seems reasonable to introduce a proper definition of a
Dirac structure on the vector bundle ${\cal
E}^1(M)=(TM\times\R)\oplus (T^\ast M\times \R)$ (a ${\cal
E}^1(M)$-Dirac structure in our terminology). This job was done by
A. Wade in \cite{Wa}. A ${\cal E}^1(M)$-Dirac structure is a
vector sub-bundle $L$ of ${\cal E}^1(M)$ that is maximally
isotropic under the natural symmetric pairing of ${\cal E}^1(M)$
and such that the space $\Gamma (L)$ is closed under a suitable
bracket $\makebox{{\bf [}}\, ,\, \makebox{{\bf ]}}$ on $\Gamma
({\cal E}^1(M))$ (this bracket may be defined using the general
algebraic constructions of Dorfman \cite{Do}). Apart from ${\cal
E}^1(M)$-Dirac structures which come from Dirac structures on $M$
or from Jacobi structures on $M$, other interesting examples can
be obtained from a homogeneous Poisson structure on $M$, from a
1-form on $M$ (a precontact structure in our terminology) or from
a locally conformal presymplectic (l.c.p.) structure, that is, a
pair $(\Omega ,\omega )$, where $\Omega$ is a 2-form on $M$,
$\omega$ is a closed 1-form and $d\Omega =\omega \wedge \Omega$
(see \cite{Wa}).

If $L$ is a ${\cal E}^1(M)$-Dirac structure, $\makebox{{\bf [}}\,
,\, \makebox{{\bf ]}}_L$ is the restriction to $\Gamma (L)\times
\Gamma (L)$ of the extended Courant bracket $\makebox{{\bf [}}\,
,\, \makebox{{\bf ]}}$ and $\rho _L$ is the restriction to L of
the canonical projection $\rho :{\cal E}^1(M)\to TM$, then the
triple $(L,\makebox{{\bf [}}\, ,\, \makebox{{\bf ]}}_L,\rho _L)$
is a Lie algebroid over $M$ (see \cite{Wa}). Using the same
terminology as in the Jacobi case, the Lie algebroid foliation
${\cal F}_L$ on $M$ associated with $L$ is called the
characteristic foliation of $L$. An important remark is that the
section $\phi _L$ of the dual bundle $L^\ast$ defined by $\phi
_L(e)=f$, for $e=(X,f)+(\alpha ,g)\in \Gamma (L)$, is a 1-cocycle
of the Lie algebroid $(L,\makebox{{\bf [}}\, ,\, \makebox{{\bf
]}}_L,\rho _L)$. Anyway, since ${\cal E}^1(M)$-Dirac structures
are closely related with Jacobi structures, it is not very
surprising the presence of a Lie algebroid and a 1-cocycle in the
theory (see \cite{IM,IM2,IM3}).

Several aspects related with the geometry of ${\cal E}^1(M)$-Dirac
structures were investigated by Wade in \cite{Wa}. However, the
nature of the induced structure on the leaves of the
characteristic foliation of a ${\cal E}^1(M)$-Dirac structure $L$
was not discussed in \cite{Wa}. So, the aim of our paper is to
describe such a nature. In addition, we will show that one may
obtain, from $L$, a Dirac structure $\tilde{L}$ on $M\times\R$ in
the sense of Courant-Weinstein and we will discuss the relation
between the induced structures on the leaves of the characteristic
foliations of $L$ and $\tilde{L}$. For the above purposes, we will
prove some general results about the relation between the
1-cocycles of an arbitrary Lie algebroid $A$ over $M$ and the
leaves of the Lie algebroid foliation on $M$ associated with $A$.
In our opinion, these last results could be of independent
interest.

The paper is organized as follows. In Section 2, we recall several
definitions and results about ${\cal E}^1(M)$-Dirac structures
which will be used in the sequel. We also present some examples
that were obtained in \cite{Wa}. In Section 3, we prove that if
$(A,\lcf \, ,\,\rcf ,\rho )$ is a Lie algebroid over $M$, $\phi$
is a 1-cocycle of $A$ and $F$ is a leaf of the Lie algebroid
foliation ${\cal F}_A$ on $M$ then $S^\phi _F=\emptyset$ or
$S^\phi _F=F$, where $S^\phi _F$ is the subset of $F$ defined by
$S^\phi _F=$ $\{ x\in F\, /\, ker\, (\rho _{|A_x})\subseteq <\phi
(x)>^\circ \}$ (see Theorem \ref{posibilidades}). Here, $A_x$ is
the fiber of $A$ over $x$ and $<\phi (x)>^\circ$ is the
annihilator of the subspace of $A^\ast _x$ generated by $\phi
(x)$. On the other hand, the Lie algebroid structure $(\lcf \,
,\,\rcf ,\rho )$ and the 1-cocycle $\phi$ induce a Lie algebroid
structure $(\lcf \,,\, \rcf \, \bar{ }\, ^{\phi},
\bar{\rho}^{\phi})$ on the vector bundle $\bar{A}=A\times \R\to
M\times\R$ (see (\ref{corchbarra})). Then, if $F$ and $\bar{F}$
are the leaves of the Lie algebroid foliations ${\cal F}_A$ and
${\cal F}_{\bar{A}}$ passing through $x_0\in M$ and $(x_0,t_0)\in
M\times\R$, we obtain, in the two possible cases ($S^\phi
_F=\emptyset$ or $S^\phi _F=F$), the relation between $F$ and
$\bar{F}$ (see Theorem \ref{rel-hojas}). Now, assume that $L$ is a
${\cal E}^1(M)$-Dirac structure and that $F$ is a leaf of the
characteristic foliation ${\cal F}_L$. Then, in Section 4, we
prove that $L$ induces, in a natural way, a ${\cal E}^1(F)$-Dirac
structure $L_F$ and, moreover (using the results of Section 3), we
describe the nature of $L_F$. In fact, we obtain that in the case
when $S^{\phi _L}_F=\emptyset$, $L_F$ comes from a precontact
structure on $F$ and in the case when $S^{\phi _L}_F=F$, $L_F$
comes from a l.c.p. structure on $F$ (see Theorem \ref{hojas}).
Using this Theorem, we directly deduce the results of Courant
\cite{Co} about the leaves of the characteristic foliation of a
Dirac structure and the results of Kirillov \cite{Ki} and
Guedira-Lichnerowicz \cite{GL} about the leaves of the
characteristic foliation of a Jacobi structure. We also apply the
theorem to the particular case when $L$ comes from a homogeneous
Poisson structure and some interesting consequences are derived.
Finally, in Section 5, we prove that a Dirac structure $\tilde{L}$
on $M\times\R$ can be obtained from a ${\cal E}^1(M)$-Dirac
structure $L$ in such a way that the Lie algebroid associated with
$\tilde{L}$ is isomorphic to the Lie algebroid over $M\times\R$,
$(\bar{L}=L\times \R,\makebox{{\bf [}}\, ,\, \makebox{{\bf ]}}_L
\kern-3pt\bar{ } \,^{\phi _L} , \bar{\rho}_L^{\phi _L})$. Thus, if
$(x_0,t_0)$ is a point of $M\times \R$, one may consider the
leaves $F$ and $\tilde{F}$ of the characteristic foliations of $L$
and $\tilde{L}$ passing through $x_0$ and $(x_0,t_0)$. Then, using
the results of Section 3, we obtain the relation between $F$ and
$\tilde{F}$ and, in addition, we describe the presymplectic 2-form
on $\tilde{F}$ in terms of the precontact structure on $F$, when
$S^{\phi _L} _F=\emptyset$, or in terms of the l.c.p. structure on
$F$, when $S^{\phi _L}_F=F$ (see Theorem \ref{relacion}). As an
application, we directly deduce some results of
Guedira-Lichnerowicz \cite{GL} about the relation between the
leaves of the characteristic foliation of a Jacobi structure on
$M$ and the leaves of the symplectic foliation of the Poisson
structure on $M\times\R$ induced by the Jacobi structure.
\section{${\cal E}^1(M)$-Dirac structures}
\setcounter{equation}{0} All the manifolds considered in this
paper are assumed to be connected and of class $C^\infty$.
Moreover, if $M$ is a differentiable manifold, we will denote by
${\cal E}^1(M)$ the vector bundle $(TM\times \R)\oplus(T^\ast
M\times \R)\to M$. Note that the space of global sections $\Gamma
({\cal E}^1(M))$ of ${\cal E}^1(M)$ can be identified with the
direct sum $(\frak X (M)\times C^\infty (M,\R ))\oplus (\Omega
^1(M) \times C^\infty (M,\R))$.
\subsection{Definition and characterization of ${\cal E}^1(M)$-Dirac
structures}\label{defi} Along this Section, we will recall the
definition of a ${\cal E}^1(M)$-Dirac structure, which was
introduced by A. Wade in \cite{Wa}. We will also give several
results related to this notion.

The natural symmetric and skew-symmetric pairings $<\, ,\, >_+$
and $<\, ,\, >_-$ on $V\oplus V^\ast$, $V$ being a real vector
space of finite dimension, can be extended, in a natural way, to
the Whitney sum $A\oplus A^\ast$, where $A\to M$ is a real vector
bundle over a manifold $M$. We also denote by $<\, ,\, >_+$ and
$<\, ,\, >_-$ the resultant pairings on $\Gamma (A\oplus
A^\ast)\cong \Gamma (A)\oplus \Gamma (A^\ast )$. In the particular
case when $A=TM\times \R$, the explicit expressions of  $<\, ,\,
>_+$ and $<\, ,\, >_-$ on $\Gamma ({\cal E}^1(M))$ are
\begin{equation}\label{emparejamientos}
\begin{array}{c}
<(X_1,f_1)+(\alpha _1,g_1),(X_2,f_2)+(\alpha
_2,g_2)>_+=\frac{1}{2}\Big ( i_{X_2}\alpha _1+f_2g_1+i_{X_1}\alpha
_2+f_1g_2\Big ),\\ <(X_1,f_1)+(\alpha _1,g_1),(X_2,f_2)+(\alpha
_2,g_2)>_-=\frac{1}{2}\Big ( i_{X_2}\alpha _1+f_2g_1-i_{X_1}\alpha
_2-f_1g_2\Big ),
\end{array}
\end{equation}
for $(X_i,f_i)+(\alpha _i,g_i)\in \Gamma ({\cal E}^1(M))$, $i\in
\{ 1,2\}$. One may also consider the homomorphism of $C^\infty
(M,\R )$-modules $\rho :\Gamma ({\cal E}^1(M))\to \frak X (M)$
defined by
\begin{equation}\label{ancla}
\rho ((X,f)+(\alpha ,g))=X.
\end{equation}
On the other hand, in \cite{Wa} A. Wade introduced a suitable
$\R$-bilinear bracket $\makebox{{\bf [}}\, ,\, \makebox{{\bf
]}}:\Gamma ({\cal E}^1(M))\times \Gamma ({\cal E}^1(M))\to \Gamma
({\cal E}^1(M))$ on the space $\Gamma ({\cal E}^1(M))$. This
approach is based on an idea that can be found in \cite{Do}, where
the author generalizes the Courant bracket on $\Gamma (TM\oplus
T^\ast M)$ to the case of complexes over Lie algebras. The bracket
$\makebox{{\bf [}}\, ,\, \makebox{{\bf ]}}$ is given by
\begin{equation}\label{corchete}
\begin{array}{ccl}
\makebox{{\bf [}}(X_1,f_1)+(\alpha _1,g_1),(X_2,f_2)+(\alpha
_2,g_2)\makebox{{\bf ]}}&=&\Big ( [X_1,X_2],X_1(f_2)-X_2(f_1)\Big
)\\ &&\kern-80pt+\Big ({\cal L}_{X_1}\alpha _2-{\cal
L}_{X_2}\alpha _1+\frac{1}{2}d(i_{X_2}\alpha _1-i_{X_1}\alpha
_2)\\ &&\kern-85pt+f_1\alpha _2-f_2\alpha
_1+\frac{1}{2}(g_2df_1-g_1df_2-f_1dg_2+f_2dg_1),\\
&&\kern-90ptX_1(g_2)-X_2(g_1)+\frac{1}{2}(i_{X_2}\alpha
_1-i_{X_1}\alpha _2-f_2g_1+f_1g_2)\Big ),
\end{array}
\end{equation}
for $(X_i,f_i)+(\alpha _i,g_i)\in \Gamma ({\cal E}^1(M))$, $i\in
\{ 1,2\}$, where $[\, ,\, ]$ is the usual Lie bracket of vector
fields and ${\cal L}$ is the Lie derivative operator on $M$. This
bracket is skew-symmetric and, moreover, we have that
\begin{equation}\label{derivacion}
\makebox{{\bf [}}e_1,f\, e_2\makebox{{\bf ]}}=f\makebox{{\bf
[}}e_1, e_2\makebox{{\bf ]}}+\rho (e_1)(f)e_2-<e_1,e_2>_+\Big (
(0,0)+(d f, 0)\Big ),
\end{equation}
for $e_1,e_2\in \Gamma ({\cal E}^1(M))$ and $f\in C^\infty (M,\R
)$. We note that $\makebox{{\bf [}}\, ,\,\makebox{{\bf ]}}$ is
not, in general, a Lie bracket, since the Jacobi identity does not
hold (see \cite{Wa}).

Now, let $L$ be a vector sub-bundle of ${\cal E}^1(M)$ which is
isotropic under the symmetric pairing $<\, ,\, >_+$. We may
consider the map $T_L:\Gamma (L)\times \Gamma (L)\times \Gamma
(L)\to C^\infty (M,\R )$ given by
\begin{equation}\label{3-tensor}
T_L(e_1,e_2,e_3)= <\makebox{{\bf [}}e_1,e_2\makebox{{\bf
]}},e_3>_+,\mbox{ for }e_1,e_2,e_3\in \Gamma (L).
\end{equation}
If $e_i=(X_i,f_i)+(\alpha _i,g_i)$ with $i\in \{ 1,2,3\}$ then,
using (\ref{emparejamientos}), (\ref{corchete}), (\ref{3-tensor})
and the fact that $<e_i,e_j>_+=0$, for $i,j\in \{1,2,3\}$, we
deduce that
\begin{equation}\label{cuentita}
\begin{array}{ccl}
T_L(e_1,e_2,e_3)=&\displaystyle\frac{1}{2}\sum_{Cycl.(e_1,e_2,e_3)}&\Big
( i_{[X_1,X_2]}\alpha _3+g_3(X_1(f_2)-X_2(f_1))\\
&&+X_3(i_{X_2}\alpha _1+f_2g_1)+f_3(i_{X_2}\alpha _1+f_2g_1)\Big
).
\end{array}
\end{equation}
Thus, from (\ref{derivacion}), (\ref{3-tensor}) and
(\ref{cuentita}), it follows that $T_L$ is a skew-symmetric
$C^\infty (M,\R)$-trilinear map, that is, $T_L\in \Gamma (\wedge
^3 L^\ast )$. Furthermore, using Proposition 3.3 in \cite{Wa}, we
obtain that
\begin{equation}\label{pseudo-Jacobi}
\makebox{{\bf [}}\makebox{{\bf [}}e_1,e_2\makebox{{\bf
]}},e_3\makebox{{\bf ]}}+\makebox{{\bf [}}\makebox{{\bf
[}}e_2,e_3\makebox{{\bf ]}},e_1\makebox{{\bf ]}}+\makebox{{\bf
[}}\makebox{{\bf [}}e_3,e_1\makebox{{\bf ]}},e_2\makebox{{\bf
]}}=(0,0)+(dT_L(e_1,e_2,e_3),T_L(e_1,e_2,e_3)),
\end{equation}
for $e_1,e_2,e_3\in \Gamma (L)$.
\begin{definition} \cite{Wa}
A ${\cal E}^1(M)$-Dirac structure on $M$ is a sub-bundle $L$ of
${\cal E}^1(M)$ which is maximally isotropic under the symmetric
pairing $<\, ,\, >_+$ and such that $\Gamma (L)$ is closed under
$\makebox{{\bf [}}\, ,\, \makebox{{\bf ]}}$.
\end{definition}
It is clear that if $L$ is a ${\cal E}^1(M)$-Dirac structure on
$M$ then the section $T_L\in \Gamma (\wedge ^3 L^\ast )$ vanishes.
In fact, we have the following result.
\begin{proposition}\label{caracterizacion} \cite{Wa}
Let $L$ be a sub-bundle of ${\cal E}^1(M)$ which is maximally
isotropic under the symmetric pairing $<\, ,\, >_+$. Then, $L$ is
a ${\cal E}^1(M)$-Dirac structure if and only if the section
$T_L\in \Gamma (\wedge ^3L^\ast )$ given by (\ref{3-tensor})
vanishes.
\end{proposition}
From (\ref{pseudo-Jacobi}) and Proposition \ref{caracterizacion},
we conclude that the restriction of $\makebox{{\bf [}}\, ,\,
\makebox{{\bf ]}}$ to $\Gamma (L)$ satisfies the Jacobi identity.
\subsection{Lie algebroids and the characteristic foliation of a ${\cal
E}^1(M)$-Dirac structure} A {\it Lie algebroid} $A$ over a
manifold $M$ is a vector bundle $A$ over $M$ together with a Lie
algebra structure $\lcf\, ,\, \rcf$ on the space $\Gamma (A)$ and
a bundle map $\rho :A \to TM$, called the {\it anchor map}, such
that if we also denote by $\rho :\Gamma (A) \to \frak X (M)$ the
homomorphism of $C^\infty (M,\R )$-modules induced by the anchor
map then:
\begin{enumerate}
\item[{\it i)}] $\rho :(\Gamma (A),\lcf \, ,\, \rcf )\to (\frak X
(M),[\, ,\, ])$ is a Lie algebra homomorphism and
\item[{\it ii)}] for all $f\in C^\infty (M,\R)$ and for all $X ,Y \in
\Gamma (A)$, one has
        $$\lcf X ,f Y \rcf=f\lcf X,Y\rcf +(\rho  (X )(f))Y.$$
\end{enumerate}
The triple $(A,\lcf \, ,\, \rcf ,\rho )$ is called a Lie algebroid
over $M$ (see \cite{Mk}).

Let $(A,\lcf \, ,\,\rcf ,\rho )$ be a Lie algebroid over $M$. We
consider the generalized distribution ${\cal F}_A$ on $M$ whose
characteristic space at a point $x\in M$ is given by
\begin{equation}\label{distribucion}
{\cal F}_A(x)=\rho (A_x),
\end{equation}
where $A_x$ is the fiber of $A$ over $x$. Since $\Gamma (A)$ is a
finitely generated $C^\infty (M,\R )$-module (see \cite{GHV}), the
distribution ${\cal F}_A$ is also finitely generated. Moreover, it
is clear that it is involutive. Thus, ${\cal F}_A$ defines a
generalized foliation on $M$ in the sense of Sussman \cite{Su}.
${\cal F}_A$ is the {\em Lie algebroid foliation} on $M$
associated with $A$.
\begin{remark}\label{curv-rang}
{\rm If $F$ is the leaf of ${\cal F}_A$ passing through $x\in M$,
$dim\,F=r$ and $y\in M$ then $y\in F$ if and only if there exists
a continuous piecewise smooth path $\gamma :I\to M$ from $x$ to
$y$, which is tangent to ${\cal F}_A$ and such that $dim\, {\cal
F}_A(\gamma (t))=r$, for all $t\in I$ (see \cite{LM,V}). }
\end{remark}
If $(A,\lcf \, ,\,\rcf ,\rho )$ is a Lie algebroid over $M$, one
can introduce the {\em Lie algebroid cohomology complex with
trivial coefficients} (for the explicit definition of this complex
we remit to \cite{Mk}). The space of 1-cochains is $\Gamma (A^\ast
)$, where $A^\ast$ is the dual bundle to $A$, and a 1-cochain
$\phi \in \Gamma (A^\ast)$ is a 1-cocycle if and only if $\phi
\lcf X,Y \rcf =\rho (X)(\phi (Y))-\rho (Y)(\phi (X)),$ for all
$X,Y\in \Gamma (A)$.

Now, suppose that $M$ is a differentiable manifold, that
$\makebox{{\bf [}}\, ,\, \makebox{{\bf ]}}$ is the bracket on
$\Gamma ({\cal E}^1(M))$ given by (\ref{corchete}) and that $\rho
:\Gamma ({\cal E}^1(M))\to \frak X(M)$ is the homomorphism of
$C^\infty (M,\R)$-modules defined by (\ref{ancla}). Assume also
that $L$ is a ${\cal E}^1(M)$-Dirac structure and that $\rho _L$
(respectively, $\makebox{{\bf [}}\, ,\, \makebox{{\bf ]}}_L$) is
the restriction of $\rho$ (respectively, $\makebox{{\bf [}}\, ,\,
\makebox{{\bf ]}}$) to $\Gamma (L)$ (respectively, $\Gamma
(L)\times \Gamma (L)$). Then, it is clear that the triple
$(L,\makebox{{\bf [}}\, ,\, \makebox{{\bf ]}}_L,\rho _L)$ is a Lie
algebroid over $M$ (see \cite{Wa} and Section \ref{defi}). Thus,
one can consider the Lie algebroid foliation ${\cal F}_L$ on $M$
associated with $L$. ${\cal F}_L$ is called the {\em
characteristic foliation of the ${\cal E}^1(M)$-Dirac structure}.

On the other hand, we may introduce a section $\phi _L$ of the
dual bundle $L^\ast$ as follows:
\begin{equation}\label{1-cociclo}
\phi _L(e)=f,\mbox{ for } e=(X,f)+(\alpha ,g)\in \Gamma (L).
\end{equation}
A direct computation, using (\ref{ancla}), (\ref{corchete}) and
(\ref{1-cociclo}), proves that $\phi _L$ is a 1-cocycle.
\subsection{Examples of ${\cal E}^1(M)$-Dirac
structures}\label{ejemplos} Next, we will present some examples of
${\cal E}^1(M)$-Dirac structures which were obtained in \cite{Wa}.
In addition, we will describe the Lie algebroids, the
characteristic foliations and the 1-cocycles associated with these
structures.

{\bf 1.- Dirac structures}

Let $\tilde{L}$ be a vector sub-bundle of $TM\oplus T^\ast M$ and
consider the vector sub-bundle $L$ of ${\cal E}^1(M)$ whose
sections are
\begin{equation}\label{ej1}
\Gamma (L)=\{ (X,0)+(\alpha ,f )/ X+\alpha \in \Gamma (\tilde{L})
,f \in C^\infty (M,\R ) \}.
\end{equation}
Then, $\tilde{L}$ {\it is a Dirac structure on }$M$ in the sense
of Courant-Weinstein \cite{Co,CW} if and only if $L$ is a ${\cal
E}^1(M)$-Dirac structure (see \cite{Wa}). We recall that a vector
sub-bundle $\tilde{L}$ of $TM\oplus T^\ast M$ is a Dirac structure
on $M$ if $\tilde{L}$ is maximally isotropic under the natural
symmetric pairing $<\, ,\, >_+$ on $TM\oplus T^\ast M$ and, in
addition, the space of sections of $L$, $\Gamma (L)$, is closed
under the {\it Courant bracket} $\makebox{{\bf [}}\, ,\,
\makebox{{\bf ]}}\,^{\tilde{ }}$ which is defined by
\begin{equation}\label{corchete-chico}
\begin{array}{ccl}
\makebox{{\bf [}}X_1+\alpha _1,X_2+\alpha _2\makebox{{\bf
]}}\,^{\tilde{ }}&=& [X_1,X_2]+ ({\cal L}_{X_1}\alpha _2-{\cal
L}_{X_2}\alpha _1+\frac{1}{2}d(i_{X_2}\alpha _1-i_{X_1}\alpha
_2)),
\end{array}
\end{equation}
for $X_1+\alpha _1, X_2+\alpha _2\in \frak X(M)\oplus \Omega
^1(M)\cong \Gamma (TM\oplus T^\ast M)$.

If $\tilde{L}\subseteq TM\oplus T^\ast M$ is a Dirac structure on
$M$ then the triple $(\tilde{L},\makebox{{\bf [}}\, ,\,
\makebox{{\bf ]}}\,^{\tilde{ }}\kern-3pt_{\tilde{L}},\tilde{\rho}
_{\tilde{L}} )$ is a Lie algebroid over $M$, where $\makebox{{\bf
[}}\, ,\, \makebox{{\bf ]}}\,^{\tilde{ }}\kern-3pt_{\tilde{L}}$ is
the restriction to $\Gamma (\tilde{L})\times \Gamma (\tilde{L})$
of the Courant bracket given by (\ref{corchete-chico}) and
$\tilde{\rho} _{\tilde{L}}$ is the restriction to $\Gamma
(\tilde{L})$ of the map $\tilde{\rho} :\Gamma (TM\oplus T^\ast
M)\to \frak X(M)$ defined by
\begin{equation}\label{ancla-chica}
\tilde{\rho} (X+\alpha )=X,
\end{equation}
for all $X+\alpha \in \Gamma (TM\oplus T^\ast M)$ (see \cite{Co}).
The {\em characteristic foliation associated with $\tilde{L}$} is
the Lie algebroid foliation ${\cal F}_{\tilde{L}}$. It is clear
that ${\cal F}_{\tilde{L}}(x)={\cal F}_L(x)$, for all $x\in M$. In
addition, from (\ref{1-cociclo}) and (\ref{ej1}), it follows that
the 1-cocycle $\phi _L$ identically vanishes.

{\bf 2.- Locally conformal presymplectic structures}

A {\it locally conformal presymplectic (l.c.p.) structure} on a
manifold $M$ is a pair $(\Omega ,\omega )$, where $\Omega$ is a
2-form on $M$, $\omega$ is a closed 1-form and $d\Omega
=\omega\wedge\Omega$. If $(\Omega ,\omega )$ is a l.c.p. structure
on $M$, one may define the vector sub-bundle $L_{(\Omega ,\omega
)}$ of ${\cal E}^1(M)$ whose sections are
\begin{equation}\label{ej3}
\Gamma (L_{(\Omega ,\omega )})=\{ (X,-i_X\omega )+(i_X\Omega
+f\,\omega, f )/(X,f )\in  \frak X(M)\times C^\infty (M,\R) \}.
\end{equation}
It is clear that the vector bundles $L_{(\Omega ,\omega )}$ and
$TM\times \R$ are isomorphic. In addition, $L_{(\Omega ,\omega )}$
is a ${\cal E}^1(M)$-Dirac structure \cite{Wa}. Note that if
$\omega =0$ then $\Omega$ is a {\it presymplectic form} on $M$.
Furthermore, if $(\Omega, \omega )$ is a l.c.p. structure on a
manifold $M$ of even dimension and $\Omega$ is a nondegenerate
2-form then $(\Omega ,\omega )$ is a {\em locally conformal
symplectic structure} (see \cite{GL,Ki,V0}).

Let $(\Omega ,\omega )$ be a l.c.p. structure on a manifold $M$
and $L_{(\Omega ,\omega )}$ be the associated ${\cal
E}^1(M)$-Dirac structure. Then, using (\ref{ancla}),
(\ref{corchete}) and (\ref{ej3}), we deduce that the Lie
algebroids $(L_{(\Omega ,\omega )},\makebox{{\bf [}}\, ,\,
\makebox{{\bf ]}}_{L_{(\Omega ,\omega )}},\rho _{L_{(\Omega
,\omega )}})$ and $(TM\times \R,\lcf \, ,\, \rcf _{(\Omega ,\omega
)} ,\rho _{(\Omega ,\omega )} )$ are isomorphic, where $\lcf \,
,\, \rcf _{(\Omega ,\omega )}$ and $\rho _{(\Omega ,\omega )}$ are
given by $$\begin{array}{l} \lcf (X,f) ,(Y,g) \rcf _{(\Omega
,\omega )}=([X,Y],\Omega (X,Y)+(X(g)-g\omega (X))-(Y(f)-f\omega
(Y))),\\
\\
\rho _{(\Omega ,\omega )}(X,f)=X,
\end{array}$$
for $(X,f),(Y,g)\in \frak X(M)\times C^\infty (M,\R)$. We remark
that the map $\nabla :\frak X(M)\times C^\infty (M,\R)\to C^\infty
(M,\R)$ defined by $\nabla _Xf=X(f)-f\omega (X)$, for $X\in \frak
X(M)$ and $f\in C^\infty (M,\R)$, induces a representation of the
Lie algebroid $(TM, [\, ,\, ],Id)$ on the trivial vector bundle
$M\times\R\to M$ and that $\Omega$ is a 2-cocycle of $(TM, [\, ,\,
],Id)$ with respect to this representation. In addition, the Lie
algebroid $(TM\times \R,\lcf \, ,\, \rcf _{(\Omega ,\omega )}
,\rho _{(\Omega ,\omega )} )$ is the extension of $(TM, [\, ,\,
],Id)$ via $\nabla$ and $\Omega$ (for the definition of the
extension of a Lie algebroid $A$ with respect to a 2-cocycle and a
representation of $A$ on a vector bundle, see \cite{Mk}). On the
other hand, it is clear that ${\cal F}_{L_{(\Omega
,\omega)}}(x)=T_xM$, for all $x\in M$, and that, under the
isomorphism between $L_{(\Omega ,\omega )}$ and $TM\times \R$, the
1-cocycle $\phi _{L_{(\Omega ,\omega )}}$ is the pair $(-\omega
,0)\in \Omega ^1(M)\times C^\infty (M,\R)\cong \Gamma (T^\ast
M\times \R)$ (see (\ref{1-cociclo}) and (\ref{ej3})).

{\bf 3.- Precontact structures}

A {\it precontact structure} on a manifold $M$ is a 1-form $\eta$
on $M$. A precontact structure $\eta$ on $M$ induces a ${\cal
E}^1(M)$-Dirac structure $L_\eta$. More precisely, suppose that
$\Phi$ is a 2-form on $M$, that $\eta$ is a 1-form and consider
the vector sub-bundle $L$ of ${\cal E}^1(M)$ whose sections are
\begin{equation}\label{ej4'}
\Gamma (L)=\{ (X,f )+(i_X\Phi +f\,\eta ,-i_X\eta )/(X,f )\in \frak
X(M)\times C^\infty (M,\R) \}.
\end{equation}
The vector bundles $L$ and $TM\times \R$ are isomorphic. Moreover,
$L$ is a ${\cal E}^1(M)$-Dirac structure if and only if $\Phi
=d\eta$ (see \cite{Wa}). Thus,
\begin{equation}\label{ej4}
\Gamma (L_\eta)=\{ (X,f )+(i_Xd\eta +f\,\eta ,-i_X\eta )/(X,f )\in
\frak X(M)\times C^\infty (M,\R) \}.
\end{equation}
Note that a precontact structure $\eta$ on a manifold $M$ of odd
dimension $2n+1$ such that $\eta \wedge (d\eta )^n$ is a volume
form is a {\it contact structure} (see \cite{GL,Ki,LM,Li2}).

Let $\eta$ be a precontact structure on a manifold $M$ and
$L_\eta$ be the associated ${\cal E}^1(M)$-Dirac structure. Then,
the Lie algebroids $(L_\eta,\makebox{{\bf [}}\, ,\, \makebox{{\bf
]}}_{L_\eta},\rho _{L_\eta})$ and $(TM\times \R ,\makebox{{\bf
[}\, ,\, {\bf ]}}, \pi)$ are isomorphic, where $\pi :TM\times \R
\to TM$ is the canonical projection over the first factor and
$\makebox{{\bf [}\, ,\, {\bf ]}}$ is the usual Lie bracket on
$\Gamma (TM\times \R)\cong \frak X(M)\times C^\infty (M,\R)$ given
by $$ \makebox{{\bf [}} (X,f),(Y,g)\makebox{{\bf
]}}=([X,Y],X(g)-Y(f)), $$ for $(X,f),(Y,g )\in \frak X (M)\times
C^\infty (M,\R )$. We also have that ${\cal F}_{L_\eta}(x)=T_xM$,
for all $x\in M$. Moreover, under the isomorphism between $L_\eta$
and $TM\times\R$, the 1-cocycle $\phi _{L_\eta}$ is the pair
$(0,1)\in \Omega^1(M)\times C^\infty (M,\R)\cong \Gamma (T^\ast
M\times \R)$ (see (\ref{1-cociclo}) and (\ref{ej4})).

{\bf 4.- Jacobi structures}

A {\it Jacobi structure} on a manifold $M$ is a pair $(\Lambda
,E)$, where $\Lambda$ is a 2-vector and $E$ is a vector field,
such that $[\Lambda ,\Lambda ]=2E\wedge \Lambda$ and $[E,\Lambda
]=0$, $[\, ,\, ]$ being the Schouten-Nijenhuis bracket. If the
vector field $E$ identically vanishes then $(M,\Lambda )$ is a
{\em Poisson manifold}. Jacobi and Poisson manifolds were
introduced by Lichnerowicz \cite{Li1,Li2} (see also
\cite{BV,DLM,GL,Ki,LM,V,We}).

Now, given a 2-vector $\Lambda$ and a vector field $E$ on a
manifold $M$, we can consider the vector sub-bundle $L_{(\Lambda
,E)}$ of ${\cal E}^1(M)$ whose sections are
\begin{equation}\label{ej2}
\Gamma (L_{(\Lambda ,E)})=\{ (\#_\Lambda (\alpha )+f\,
E,-i_E\alpha )+(\alpha ,f )/ (\alpha,f ) \in \Omega ^1(M)\times
C^\infty (M,\R) \},
\end{equation}
where $\# _\Lambda :\Omega ^1(M)\to \frak X (M)$ is the
homomorphism of $C^\infty (M,\R)$-modules defined by $\beta (\#
_\Lambda (\alpha ))=\Lambda (\alpha ,\beta )$, for $\alpha ,\beta
\in \Omega ^1(M)$. Note that the vector bundles $L_{(\Lambda ,E)}$
and $T^\ast M\times \R$ are isomorphic. Moreover, we have that
$L_{(\Lambda ,E)}$ is a ${\cal E}^1(M)$-Dirac structure if and
only if $(\Lambda ,E)$ is a Jacobi structure (see \cite{Wa}).

If $(\Lambda ,E)$ is a Jacobi structure on a manifold $M$ and
$L_{(\Lambda ,E)}$ is the associated ${\cal E}^1 (M)$-Dirac
structure then the Lie algebroids $(L_{(\Lambda ,E)},\makebox{{\bf
[}}\, ,\, \makebox{{\bf ]}}_{L_{(\Lambda ,E)}},\rho _{L_{(\Lambda
,E)}})$ and $(T^\ast M \times \R ,$ $\lcf \, ,\, \rcf _{(\Lambda
,E)},$ $\widetilde{\#}_{(\Lambda ,E)})$ are isomorphic, where
$\lcf \, ,\,\rcf _{(\Lambda ,E)}$ and
$\widetilde{\#}_{(\Lambda,E)}$ are defined by $$
\begin{array}{cll}
\kern-17.6pt\lcf (\alpha ,f),(\beta ,g)\rcf _{(\Lambda
,E)}&\kern-9pt=&\kern-9pt({\cal L}_{\#_{\Lambda}(\alpha
)}\beta\kern-2pt-\kern-2pt{\cal L}_{\#_{\Lambda}(\beta
)}\alpha\kern-2pt -\kern-2pt d(\Lambda (\alpha ,\beta
))\kern-2pt+\kern-2ptf{\cal L}_{E}\beta\kern-2pt-\kern-2ptg {\cal
L}_{E}\alpha \kern-2pt-\kern-2pti_{E}(\alpha \wedge \beta),\\
&\kern-9pt &\kern-9pt \Lambda (\beta ,\alpha
)\kern-2pt+\kern-2pt\#_{\Lambda}(\alpha
)(g)\kern-2pt-\kern-2pt\#_{\Lambda}(\beta
)(f)\kern-2pt+\kern-2ptfE(g)\kern-2pt-\kern-2ptg E(f)),\\ & & \\
\kern-17.6pt\widetilde{\#}_{(\Lambda ,E)}(\alpha
,f)&\kern-9pt=&\#_{\Lambda}(\alpha )+f E,
\end{array}
$$ for $(\alpha ,f),(\beta ,g)\in \Omega ^1(M)\times C^\infty
(M,\R )\cong \Gamma (T^\ast M\times \R)$ (see \cite{Wa}). The Lie
algebroid structure $(\lcf \, ,\,\rcf _{(\Lambda ,E)},
\widetilde{\#}_{(\Lambda,E)})$ on $T^\ast M\times \R$ was
introduced in \cite{KS}. In this case, the characte\-ristic
foliation of the ${\cal E}^1(M)$-Dirac structure $L_{(\Lambda
,E)}$ is just the characteristic foliation on $M$ associated with
the Jacobi structure $(\Lambda ,E)$ (see \cite{DLM,GL,Ki}). In
addition, under the isomorphism between $L_{(\Lambda ,E)}$ and
$T^\ast M\times \R$, the 1-cocycle $\phi _{L_{(\Lambda ,E)}}$ is
the pair $(-E,0)\in \frak X(M)\times C^\infty (M,\R)\cong \Gamma
(TM\times \R)$ (see (\ref{1-cociclo}) and (\ref{ej2})).

{\bf 5.- Homogeneous Poisson structures}

A {\it homogeneous Poisson manifold} $(M,\Pi, Z)$ is a Poisson
manifold $(M,\Pi )$ with a vector field $Z$ such that $[Z,\Pi
]=-\Pi$ (see \cite {DLM}). Given a 2-vector $\Pi$ and a vector
field $Z$ on a manifold $M$, we can define the vector sub-bundle
$L_{(\Pi ,Z)}$ of ${\cal E}^1(M)$ whose sections are
\begin{equation}\label{ej5}
\Gamma (L_{(\Pi ,Z)})=\{ (\#_ \Pi (\alpha )-f\,Z,f )+(\alpha
,i_Z\alpha )/(\alpha ,f )\in \Omega ^1(M)\times C^\infty (M,\R)\}.
\end{equation}
The vector bundles $L_{(\Pi ,Z)}$ and $T^\ast M\times\R$ are
isomorphic. Furthermore, $(M,\Pi ,Z)$ is a homogeneous Poisson
manifold if and only if $L_{(\Pi ,Z)}$ is a ${\cal E}^1(M)$-Dirac
structure (see \cite{Wa}).

Let $(M,\Pi ,Z)$ be a homogeneous Poisson manifold and $L_{(\Pi
,Z)}$ be the associated ${\cal E}^1(M)$-Dirac structure. Then,
from (\ref{ancla}), (\ref{corchete}) and (\ref{ej5}), it follows
that the Lie algebroids $(L_{(\Pi ,Z)},\makebox{{\bf [}}\, ,\,
\makebox{{\bf ]}}_{L_{(\Pi ,Z)}},\rho _{L_{(\Pi ,Z)}})$ and
$(T^\ast M \times \R ,\lcf \, ,\, \rcf _{(\Pi ,Z)},$
$\widetilde{\#}_{(\Pi ,Z)})$ are isomorphic, where $\lcf \,
,\,\rcf _{(\Pi ,Z)}$ and $\widetilde{\#}_{(\Pi ,Z)}$ are defined
by
\begin{equation}\label{par-mach}
\begin{array}{cll}
\kern-14.5pt\lcf (\alpha ,f),(\beta ,g)\rcf _{(\Pi
,Z)}&\kern-9pt=&\kern-9pt\Big ({\cal L}_{\#_{\Pi}(\alpha
)}\beta\kern-2pt-\kern-2pt{\cal L}_{\#_{\Pi}(\beta
)}\alpha\kern-2pt -\kern-2pt d(\Pi (\alpha ,\beta
))\kern-2pt-\kern-2pt f({\cal L}_{Z}\beta-\beta)+g ({\cal
L}_{Z}\alpha -\alpha),\\ &\kern-9pt &\kern-9pt \#_{\Pi}(\alpha
)(g)\kern-2pt-\kern-2pt\#_{\Pi}(\beta
)(f)\kern-2pt+g\,Z(f)\kern-2pt-\kern-2pt f\,Z(g)\Big ),\\ & &
\\ \kern-15pt\widetilde{\#}_{(\Pi ,Z)}(\alpha
,f)&\kern-9pt=&\#_{\Pi}(\alpha )-f\,Z,
\end{array}
\end{equation}
for $(\alpha ,f),(\beta ,g)\in \Omega ^1(M)\times C^\infty (M,\R
)\cong \Gamma (T^\ast M\times \R)$. Furthermore, it is clear that
${\cal F}_{L_{(\Pi ,Z)}}(x)=\# _\Pi (T^\ast _xM)+<Z(x)>$, for all
$x\in M$. In other words, if ${\cal F}_\Pi$ is the symplectic
foliation of the Poisson manifold $(M,\Pi )$ then
\begin{equation}\label{symp+Z}
{\cal F}_{L_{(\Pi ,Z)}}(x)={\cal F}_\Pi(x)+<Z(x)>,\mbox{ for all
}x\in M.
\end{equation}
In addition, under the isomorphism between $L_{(\Pi ,Z)}$ and
$T^\ast M\times\R$, the 1-cocycle $\phi_{L_{(\Pi ,Z)}}$ is the
pair $(0,1)\in \frak X(M)\times C^\infty (M,\R)\cong \Gamma
(TM\times \R)$ (see (\ref{1-cociclo}) and (\ref{ej5})).

On the other hand, we may consider the Lie algebroid structure
$(\lcf \, ,\, \rcf _\Pi ,\# _\Pi )$ on the vector bundle $T^\ast
M\to M$ induced by the Poisson structure $\Pi$ and the Lie
algebroid structure $(\lcf \, ,\, \rcf _Z,\rho _Z)$ on the vector
bundle $M\times \R\to M$ induced by the vector field $Z$. The
explicit definitions of $\lcf \, ,\, \rcf _\Pi$, $\lcf \, ,\, \rcf
_Z$ and $\rho _Z$ are $$\begin{array}{l}\kern-60pt \lcf \alpha
,\beta \rcf _\Pi ={\cal L}_{\#_\Pi (\alpha )}\beta -{\cal
L}_{\#_\Pi (\beta )}\alpha -d(\Pi (\alpha ,\beta)),\\
\kern-60pt\lcf f ,g \rcf _Z =g\,Z(f)-f\,Z(g),\quad \rho
_Z(f)=-f\,Z,
\end{array}$$
for $\alpha ,\beta \in \Omega ^1(M)$ and $f,g\in C^\infty (M,\R)$.
Then, using (\ref{par-mach}), we conclude that $(T^\ast M,\lcf \,
,\, \rcf _\Pi ,$ $\# _\Pi )$ and $(M\times \R,\lcf \, ,\, \rcf
_Z,\rho _Z)$ are Lie subalgebroids of $(T^\ast M \times \R ,\lcf
\, ,\, \rcf _{(\Pi ,Z)},$ $\widetilde{\#}_{(\Pi ,Z)})$. This
implies that $(T^\ast M,\lcf \, ,\, \rcf _\Pi ,\# _\Pi )$ and
$(M\times \R,\lcf \, ,\, \rcf _Z,\rho _Z)$ form a matched pair of
Lie algebroids in the sense of Mokri \cite{Mo}.

In Section \ref{principal}, we will prove that if $F$ is a leaf of
the characteristic foliation of a ${\cal E}^1(M)$-Dirac structure
$L$ then $L$ induces a ${\cal E}^1(F)$-Dirac structure $L_F$ and,
in addition, we will describe the nature of $L_F$. First, in the
next Section, we will show two general results about the relation
between the 1-cocycles of an arbitrary Lie algebroid $A$ and the
leaves of the Lie algebroid foliation ${\cal F}_A$.

\section{1-cocycles of a Lie algebroid and the leaves of the Lie
algebroid foliation} \setcounter{equation}{0} Let $(A,\lcf \, , \,
\rcf, \rho )$ be a Lie algebroid over $M$, $\phi\in \Gamma (A^\ast
)$ be a 1-cocycle and $\pi _1 :M\times \R \to M$ be the canonical
projection onto the first factor. We consider the map $\cdot
:\Gamma (A)\times C^\infty (M\times \R ,$ $\R )\to C^\infty
(M\times \R,\R )$ given by $$X\cdot \bar{f}=\rho(X)(\bar{f})+\phi
(X)\frac{\partial \bar{f}}{\partial t}.$$ It is easy to prove that
$\cdot$ is an action of $A$ on $M\times \R$ in the sense of
\cite{HM} (see Definition 2.3 in \cite{HM}). Thus, if $\pi _1^\ast
A$ is the pull-back of $A$ over $\pi _1$ then the vector bundle
$\pi _1^\ast A\to M\times \R$ admits a Lie algebroid structure
$(\lcf \,,\, \rcf \, \bar{ }\, ^\phi , \bar{\rho}^\phi)$ (see
Theorem 2.4 in \cite{HM}). For the sake of simplicity, when the
1-cocycle $\phi$ is zero, we will denote by $(\lcf \,,\, \rcf \,
\bar{ }, \bar{\rho})$ the resultant Lie algebroid structure on
$\pi _1^\ast A\to M\times \R$. On the other hand, it is clear that
the vector  bundles $\pi _1^\ast A\to M\times \R$ and $\bar{A}=
A\times \R\to M\times \R$ are isomorphic and that the space of
sections $\Gamma (\bar{A})$ of $\bar{A}\to M\times \R$ can be
identified with the set of time-dependent sections of $A\to M$.
Under this identification, we have that $\lcf \bar{X},\bar{Y} \rcf
\, \bar{ }\,(x,t)=\lcf \bar{X}_t,\bar{Y}_t\rcf (x)$ and that
$\bar{\rho}(\bar{X})(x,t)=\rho (\bar{X}_t)(x)$, for
$\bar{X},\bar{Y}\in \Gamma (\bar{A})$ and $(x,t)\in M\times \R$
(see \cite{HM}). In addition,
\begin{equation}\label{corchbarra}
\lcf \bar{X},\bar{Y}\rcf \, \bar{ }\, ^\phi = \lcf \bar{X},
\bar{Y}\rcf\,\bar{ } + \phi (\bar{X})\frac{\partial
\bar{Y}}{\partial t}-\phi (\bar{Y})\frac{\partial
\bar{X}}{\partial t},\qquad \bar{\rho}^\phi (\bar{X})=
\bar{\rho}(\bar{X})+\phi (\bar{X}) \frac{\partial}{\partial t},
\end{equation}
where $\displaystyle \frac{\partial\bar{X}}{\partial t}\in \Gamma
(\bar{A})$ denotes the derivative of $\bar{X}$ with respect to the
time.

Now, if ${\cal F}_{\bar{A}}$ is the Lie algebroid foliation of
$(\bar{A},\lcf \,,\, \rcf \, \bar{ }\, ^\phi , \bar{\rho}^\phi )$
then, from (\ref{corchbarra}), it follows that
\begin{equation}\label{folbarra}
{\cal F}_{\bar{A}}(x,t)=\{ \rho (e_x)+\phi
(x)(e_x)\frac{\partial}{\partial t}_{|t}\in T_{(x,t)}(M\times
\R)\, /\, e_x\in A_x\} ,
\end{equation}
for all $(x,t)\in M\times \R$. Moreover, a direct computation
shows that,
\begin{equation}\label{desigualdad}
dim\,{\cal F}_A(x)\leq dim\,{\cal F}_{\bar{A}}(x,t)\leq dim\,{\cal
F}_A(x)+1,
\end{equation}
\begin{equation}\label{equivalencia}
dim\, {\cal F}_{\bar{A}}(x,t)=dim\,{\cal F}_A(x)\iff ker\,(\rho
_{|A_x})\subseteq <\phi (x)>^\circ ,
\end{equation}
where ${\cal F}_A$ is the Lie algebroid foliation of $A$ and
$<\phi (x)>^\circ$ is the annihilator of the subspace of $A^\ast
_x$ generated by $\phi (x)$, that is, $$<\phi (x)>^\circ=\{ e_x\in
A_x\, /\, \phi (x)(e_x)=0\}.$$

\vspace{-.4cm}
\begin{remark}
{\rm Note that the vector field $\frac{\partial}{\partial t}$ on
$M\times\R$ is an infinitesimal automorphism of the foliation
${\cal F}_{\bar A}$. Therefore, if $(x,t_0)$ and $(x,t'_0)$ are
points of $M\times \R$ and $\bar{F}$, $\bar{F}'$ are the leaves of
${\cal F}_{\bar{A}}$ passing through $(x,t_0)$ and $(x,t'_0)$,
then the map $(y,s)\mapsto (y,s+(t'_0-t_0))$ is a diffeomorphism
from $\bar{F}$ to $\bar{F}'$.}
\end{remark}
Next, we will discuss some relations between the leaves of ${\cal
F}_A$ and the 1-cocycle $\phi$ and between the leaves of ${\cal
F}_A$ and ${\cal F}_{\bar{A}}$. More precisely, the aim of this
Section is to prove the following two results.
\begin{theorem}\label{posibilidades}
Let $(A,\lcf \, ,\,\rcf ,\rho )$ be a Lie algebroid and $\phi \in
\Gamma (A^\ast )$ be a 1-cocycle. If $F$ is a leaf of the Lie
algebroid foliation ${\cal F}_A$ and $S^\phi _F$ is the subset of
$F$ defined by $$S^\phi _F=\{ x\in F\, /\, ker\, (\rho
_{|A_x})\subseteq <\phi (x)>^\circ \},$$ then $S^\phi
_F=\emptyset$ or $S^\phi _F=F$. Furthermore, in the second case
$(S^\phi _F=F)$, the 1-cocycle $\phi$ induces a closed 1-form
$\omega _F$ on $F$ which is characterized by the condition
\begin{equation}\label{1-form-induc}
\omega _F(\rho (e)_{|F})=-\phi (e)_{|F},\mbox{ for all }e\in
\Gamma (A).
\end{equation}
\end{theorem}
\begin{theorem}\label{rel-hojas}
Let $(A,\lcf \, ,\,\rcf ,\rho )$ be a Lie algebroid, $\phi \in
\Gamma (A^\ast )$ be a 1-cocycle and consider on the vector bundle
$\bar{A}=A\times\R\to M\times \R$ the Lie algebroid structure
$(\lcf \,,\, \rcf \, \bar{ }\, ^\phi , \bar{\rho}^\phi )$ given by
(\ref{corchbarra}). Suppose that $(x_0,t_0)\in M\times\R$ and that
$F$ and $\bar{F}$ are the leaves of the Lie algebroid foliations
${\cal F}_A$ and ${\cal F}_{\bar{A}}$ passing through $x_0\in M$
and $(x_0,t_0)\in M\times \R$, respectively. Then:

{\it i)} If $ker\,(\rho _{|A_{x_0}})\not\subseteq <\phi
(x_0)>^\circ$ (or, equivalently, $S^\phi _F=\emptyset$) we have
that $\bar{F}=F\times\R$.

{\it ii)} If $ker\,(\rho _{|A_{x_0}})\subseteq <\phi (x_0)>^\circ$
(or, equivalently, $S^\phi _F=F$) and $\pi _1:M\times \R\to M$ is
the canonical projection onto the first factor, we have that $\pi
_1(\bar{F})=F$ and that the map $\pi _1{ } _{|\bar{F}}: \bar{F}\to
F$ is a covering map. In addition, if $\,\omega _F$ is the closed
1-form on $F$ characterized by the condition (\ref{1-form-induc})
and $i_{\bar{F}}:\bar{F}\to M\times \R$ is the canonical inclusion
then $\bar{F}$ is diffeomorphic to a Galois covering of $F$
associated with $\omega _F$ and $$(\pi _1{}_{|\bar{F}})^\ast
(\omega _F)=-d((i_{\bar{F}})^\ast (t)).$$
\end{theorem}
In order to prove the above results, we will use the following
three lemmas.
\begin{lemma}\label{lema1}
{\it i)} If $\bar{F}$ is a leaf of the Lie algebroid foliation
${\cal F}_{\bar{A}}$ and $S^\phi _{\bar{F}}$ is the subset of
$\bar{F}$ defined by $$S^\phi _{\bar{F}}=\{ (x,t)\in \bar{F}\, /\,
ker\, (\rho _{|A_x})\subseteq <\phi (x)>^\circ \},$$ then $S^\phi
_{\bar{F}}$ is an open subset of $\bar{F}$.

{\it ii)} If $F$ is a leaf of the Lie algebroid foliation ${\cal
F}_A$ then $S^\phi _F$ is an open subset of $F$.
\end{lemma}
\prueba {\it i)} Assume that $dim\, \bar{F}=r$ and let $(x_0,t_0)$
be a point of $S^\phi _{\bar{F}}$. We will show that there exists
a connected open subset $\bar{W}_{(x_0,t_0)}$ of $\bar{F}$ such
that $(x_0,t_0)\in \bar{W}_{(x_0,t_0)}$ and
$\bar{W}_{(x_0,t_0)}\subseteq S^\phi _{\bar{F}}$ or, equivalently
(see (\ref{equivalencia})),
\begin{equation}\label{conddim}
dim\,{\cal F}_A(x)=dim\,{\cal F}_{\bar{A}}(x,t)=r,\mbox{ for all }
(x,t)\in \bar{W}_{(x_0,t_0)}.
\end{equation}
Note that
\begin{equation}\label{cond2}
dim\, {\cal F}_{\bar{A}}(x,t)=r,\mbox{ for all }(x,t)\in \bar{F}.
\end{equation}
Therefore, we can choose a global generator system $\{ e_1,\ldots
,e_m\}$ of $\Gamma (A)$ in such a way that the set of vectors $\{
\rho (e_i(x_0))+\phi (x_0)(e_i(x_0))\frac{\partial}{\partial
t}_{|t_0}\}_{1\leq i\leq r}$ are a basis of the vector space
${\cal F}_{\bar{A}}(x_0,t_0)$. Then, using that $(x_0,t_0)\in
S^\phi _{\bar{F}}$, we deduce that the vectors $\{\rho
(e_i(x_0))\}_{1\leq i\leq r}$ are linearly independent in
$T_{x_0}M$. This implies that $dim\, {\cal F}_A(x_0)\geq r$ and,
since the rank of a differentiable generalized distribution is a
lower semicontinuous function (see \cite{V}), there exists an open
subset $V'_{x_0}$ of $M$, $x_0\in V'_{x_0}$, such that
\begin{equation}\label{dimension}
dim\, {\cal F}_A(x)\geq dim\, {\cal F}_A(x_0)\geq r,\mbox{ for all
}x\in V'_{x_0}.
\end{equation}
Thus, if $\bar{W}_{(x_0,t_0)}$ is the connected component of the
point $(x_0,t_0)$ in the open subset (of $\bar{F}$) $\bar{F}\cap
(V'_{x_0}\times\R)$ then, from (\ref{desigualdad}), (\ref{cond2})
and (\ref{dimension}), it follows that (\ref{conddim}) holds.

{\it ii)} Let $x_0$ be a point of $S^\phi _F$. We will prove that
there exists an open neighbourhood $V_{x_0}$ of $x_0$ in $F$ which
is contained in $S^\phi _F$.

Suppose that $t_0\in \R$, that $\bar{F}$ is the leaf of ${\cal
F}_{\bar{A}}$ passing through $(x_0,t_0)$ and that the dimension
of $\bar{F}$ is $r$. Then, the point $(x_0,t_0)\in S^\phi
_{\bar{F}}$ and thus there exists a connected open subset
$\bar{W}_{(x_0,t_0)}$ of $\bar{F}$ such that $(x_0,t_0)\in
\bar{W}_{(x_0,t_0)}$ and (\ref{conddim}) holds.

Now, if $\pi _1:M\times \R\to M$ is the canonical projection onto
the first factor, (\ref{folbarra}), (\ref{conddim}) and Remark
\ref{curv-rang} imply that $\pi _1(\bar{W}_{(x_0,t_0)}) \subseteq
F$ and that $\pi _1:\bar{W}_{(x_0,t_0)}\to F$ is a local
diffeomorphism. Therefore, $V_{x_0}=\pi _1(\bar{W}_{(x_0,t_0)})$
is an open subset of $F$, $x_0\in V_{x_0}$ and it is clear that
$V_{x_0}\subseteq S^\phi _F$ (see (\ref{equivalencia}) and
(\ref{conddim})).\QED
\begin{lemma}\label{lema2}
If $F$ is a leaf of the Lie algebroid foliation ${\cal F}_A$ then
$S^\phi _F$ is a closed subset of $F$.
\end{lemma}
\prueba We must prove that $F-S^\phi _F$ is an open subset of $F$.

Let $x_0$ be a point of $F-S^\phi _F$.

Assume that $\{ s_1,\ldots ,s_n\}$ is a local basis of $\Gamma
(A)$ in an open neighbourhood of $x_0$ in $M$. If the dimension of
$F$ is $r$ we can suppose, without loss of the generality, that
the vectors $\{ \rho (s_i(x_0))\}_{1\leq i\leq r}$ are linearly
independent in $T_{x_0}M$. Then, there exists an open subset
$U'_{x_0}$ of $M$, $x_0\in U'_{x_0}$, such that for every $x\in
U'_{x_0}$ the set of vectors $\{ \rho (s_i(x))\}_{1\leq i\leq r}$
are linearly independent in $T_xM$. Thus, this set is a basis of
$T_xF$, for all $x\in U_{x_0}=U'_{x_0}\cap F$. This implies that
$$\rho (s_{r+k})_{|U_{x_0}} = \displaystyle\sum_{i=1}^r
g^i_k\,\rho (s_i )_{|U_{x_0}},$$ for all $k\in \{ 1,\ldots
,n-r\}$, with $g^i_k\in C^\infty (U_{x_0}, \R)$. Therefore, if
$k\in \{ 1,\ldots ,n-r\}$ and $e_{r+k}:U_{x_0}\to A$ is the map
defined by $$e_{r+k}(x)=s_{r+k}(x)-\displaystyle\sum_{i=1}^r
g^i_k(x) s_i(x),\mbox{ for all }x\in U_{x_0},$$ we deduce that the
set $\{ e_{r+k}(x)\}_{1\leq k\leq n-r}$ is a basis of the vector
space $ker\, (\rho _{|A_x})$, for all $x\in U_{x_0}$.
Consequently, since $x_0\in F-S^\phi _F$, it follows that $\phi
(x_0)(e_{r+k_0}(x_0))\neq 0$, for some $k_0 \in\{ 1,\ldots
,n-r\}$. Hence, we conclude that there exists an open subset
$V_{x_0}$ of $F$, $x_0\in V_{x_0}\subseteq U_{x_0}$, such that
$\phi (x)(e_{r+k_0}(x))\neq 0$, for all $x\in V_{x_0}$. Then, it
is clear that $V_{x_0}\subseteq F-S^\phi _F$.\QED
\begin{lemma}\label{lema3}
If $\bar{F}$ is a leaf of the Lie algebroid foliation ${\cal
F}_{\bar{A}}$ and $S^\phi _{\bar{F}}$ is the subset of $\bar{F}$
defined by $$S^\phi _{\bar{F}}=\{ (x,t)\in \bar{F}\, /\, ker\,
(\rho _{|A_x})\subseteq <\phi (x)>^\circ \},$$ then $S^\phi
_{\bar{F}}$ is a closed subset of $\bar{F}$.
\end{lemma}
\prueba Assume that the dimension of $\bar{F}$ is $r$. We will
prove that $\bar{F}-S^\phi _{\bar{F}}$ is an open subset of
$\bar{F}$.

Let $(x_0,t_0)$ be a point of $\bar{F}-S^\phi _{\bar{F}}$ and
suppose that $F$ is the leaf of ${\cal F}_A$ passing through
$x_0$. We have that $x_0\in F-S^\phi _F$. Thus, using Lemmas
\ref{lema1} and \ref{lema2} and the fact that $F$ is connected, we
deduce that
\begin{equation}\label{vacio}
S^\phi _F=\emptyset.
\end{equation}
Therefore, from (\ref{folbarra}) and (\ref{vacio}), it follows
that $${\cal F}_{\bar{A}}(x,t)={\cal F}_A(x)\oplus
<\frac{\partial}{\partial t}_{|t}>,\mbox{ for all }(x,t)\in
F\times\R.$$ Consequently, $F\times\R$ is a connected integral
submanifold of ${\cal F}_{\bar{A}}$ and its dimension is $r$. This
implies that $F\times\R$ is an open subset of ${\bar F}$. Finally,
from (\ref{vacio}), we conclude that $F\times \R\subseteq
\bar{F}-S^\phi _{\bar{F}}$.\QED

\vspace{.25cm} {\it Proof of Theorem \ref{posibilidades}:} Using
Lemmas \ref{lema1} and \ref{lema2} and the fact that $F$ is
connected it follows that $S^\phi _F=\emptyset$ or $S^\phi _F=F$.

Now, if $S^\phi _F=F$ we may introduce a 1-form $\omega _F$ on $F$
given by $$\omega _F(x)(\rho (e_x))=-\phi (x)(e_x),$$ for all
$x\in F$ and $e_x\in A_x$. Note that the condition $$ker\, (\rho
_{|A_x})\subseteq <\phi (x)>^\circ,\mbox{ for all }x\in F,$$
implies that $\omega _F(x):T_xF\to \R$ is well defined. Moreover,
it is clear that $\omega _F$ satisfies (\ref{1-form-induc}) and,
since $\phi$ is a 1-cocycle, we deduce that $\omega _F$ is
closed.\QED

{\it Proof of Theorem \ref{rel-hojas}:} {\it i)} If $ker\, (\rho
_{|A_{x_0}})\not \subseteq <\phi (x_0)>^\circ$ then $(x_0,t_0)\in
\bar{F}-S^\phi _{\bar{F}}$ and thus, using Lemmas \ref{lema1} and
\ref{lema3} and the fact that $\bar{F}$ is connected, we obtain
that
\begin{equation}\label{posib1}
S^\phi _{\bar{F}}=\emptyset.
\end{equation}
Now, proceeding as in the proof of Lemma \ref{lema3}, we deduce
that $F\times\R$ is an open subset of $\bar{F}$. On the other
hand, from (\ref{desigualdad}), (\ref{equivalencia}) and
(\ref{posib1}), it follows that $$dim\,{\cal F}_A(x)=dim\,{\cal
F}_{\bar{A}}(x,t)-1=dim\,\bar{F}-1,\mbox{ for all
}(x,t)\in\bar{F}.$$ Using this fact, (\ref{folbarra}) and Remark
\ref{curv-rang}, we have that $\pi _1(\bar{F})\subseteq F$.
Therefore, we have proved that $\bar{F}=F\times \R$.

{\it ii)} Assume that $ker\, (\rho _{|A_{x_0}})\subseteq <\phi
(x_0)>^\circ$. Then, $(x_0,t_0)\in S^\phi _{\bar{F}}$ which
implies that $\bar{F}=S^\phi _{\bar{F}}$, that is,
\begin{equation}\label{posib2}
dim\,{\cal F}_A(x)=dim\,{\cal
F}_{\bar{A}}(x,t)=dim\,\bar{F},\mbox{ for all }(x,t)\in\bar{F}.
\end{equation}
Using (\ref{folbarra}), (\ref{posib2}) and Remark \ref{curv-rang},
we obtain that $\pi _1(\bar{F})\subseteq F$ and that $\pi _1{
}_{|\bar{F}}:\bar{F}\to F$ is a local diffeomorphism.
Consequently, $\pi _1(\bar{F})$ is an open subset of $F$.

In addition, from (\ref{folbarra}) and (\ref{1-form-induc}), it
follows that
\begin{equation}\label{otrarel}
(\pi _1{ }_{|\bar{F}})^\ast(\omega _F)=-d((i_{\bar{F}})^\ast t).
\end{equation}
Next, we will show that $\pi _1(\bar{F})$ is a closed subset of
$F$ and that $\pi _1{ }_{|\bar{F}}:\bar{F}\to F$ is a covering
map.

Let $x$ be a point of $F$. Since $\omega _F$ is a closed 1-form,
there exists a connected open subset $U$ in $F$ and a real
$C^\infty$-differentiable function $f_F$ on $U$ such that $x\in U$
and
\begin{equation}\label{Poincare}
\omega _F=df_F \quad\mbox{ on }U.
\end{equation}
Then, using (\ref{folbarra}), (\ref{1-form-induc}),
(\ref{Poincare}), Remark \ref{curv-rang} and the fact that $\pi
_1(\bar{F})\subseteq F$, we deduce the following result
\begin{equation}\label{implicacion}
\begin{array}{ccl}
(y,s)\in (\pi _1{ }_{|\bar{F}})^{-1}(U)=\pi _1^{-1}(U)\cap
\bar{F}&\Longrightarrow& \\ & &\kern-70pt \{
(z,s+f_F(y)-f_F(z))\in M\times \R\, /\, z\in U\}\subseteq \bar{F}.
\end{array}
\end{equation}
Thus, if $x\in F-\pi_1(\bar{F})$, we have that $U\subseteq F-\pi
_1(\bar{F})$. This proves that $\pi _1(\bar{F})$ is a closed
subset of $F$ which implies that $\pi _1(\bar{F})=F$.

Now, suppose that $(x,t)$ is a point of $\bar{F}$ and let $U$ be a
connected open subset of $F$ and $f_F$ be a real
$C^\infty$-differentiable on $F$ such that $x\in U$ and
(\ref{Poincare}) holds. If $C_{(y,s)}$ is the connected component
of a point $(y,s)\in (\pi _1{ }_{|\bar{F}})^{-1}(U)$ then, using
(\ref{otrarel}) and (\ref{Poincare}), it follows that the function
$(\pi _1{ }_{|\bar{F}})^\ast(f_F)+(i_{\bar{F}})^\ast t$ is
constant on $C_{(y,s)}$. Therefore, from (\ref{implicacion}), we
obtain that $$C_{(y,s)}=\{ (z,s+f_F(y)-f_F(z))\in M\times \R\, /\,
z\in U\}.$$ Consequently, the map $\pi _1{
}_{|C_{(y,s)}}:C_{(y,s)}\to U$ is a diffeomorphism. This proves
that $\pi _1{ }_{|\bar{F}}:\bar{F}\to F$ is a covering map.

Finally, let $E$ be the covering of $F$ associated with $\omega
_F$, that is, $E$ is the sheaf of germs of $C^\infty$ functions
$g_F$ on $F$ such that $dg_F=\omega _F$ (see Section 2 of Chapter
XIV in \cite{Go}). Denote by $(f^0_F)_{[x_0]}$ the germ of $f^0_F$
at $x_0$, where $f^0_F$ is a $C^\infty$ function on a connected
open subset $U_0$ of $F$ such that $x_0\in U_0$,
$(f^0_F)(x_0)=t_0$ and $\omega _F{ }_{|U_0}=df^0_F$. Then, using
the above description of the leaf $\bar{F}$ and the results in
\cite{Go}, we deduce that $\bar{F}$ is diffeomorphic to the
connected component of $(f^0_F)_{[x_0]}$ in $E$. In other words,
$\bar{F}$ is diffeomorphic to a Galois covering of $F$ associated
with $\omega _F$.\QED
\section{${\cal E}^1(M)$-Dirac structures, submanifolds of the
base space and the leaves of the characteristic
foliation}\label{principal}\setcounter{equation}{0}
\subsection{${\cal E}^1(M)$-Dirac structures and submanifolds of the
ba\-se space}\label{induc-base} In this Section, we will prove
that if $S$ is a submanifold of $M$ then, under certain regularity
conditions, a ${\cal E}^1(M)$-Dirac structure induces a ${\cal
E}^1(S)$-Dirac structure. This result will be used in Section
\ref{subprincipal}.

Let $L$ be a vector sub-bundle of ${\cal E}^1(M)$ which is
maximally isotropic under the symmetric pairing $<\, ,\, >_+$ and
$S$ be a submanifold of $M$. If $x$ is a point of $S$, we may
define the vector space $(L_S)_x$ by
\begin{equation}\label{subespacio}
(L_S)_x=\frac{L_x\cap ((T_xS\times\R)\oplus (T^\ast
_xM\times\R))}{L_x\cap (\{ 0\}\oplus ((T_xS)^\circ \times\{
0\}))},
\end{equation}
where $(T_xS)^\circ$ is the annihilator of $T_xS$, that is,
$(T_xS)^\circ=\{ \alpha \in T^\ast _xM\, /\, \alpha _{|T_xS}=0\}$.
We have that the linear map $(L_S)_x\to (T_xS\times \R)\oplus
(T^\ast _xS\times \R)$ given by
\begin{equation}\label{inclusion}
[(u,\lambda )+(\alpha ,\mu )]\mapsto (u,\lambda )+(\alpha
_{|T_xS},\mu ),
\end{equation}
is a monomorphism and thus $(L_S)_x$ can be identified with a
subspace of $(T_xS\times \R)\oplus (T^\ast _xS\times\R )$.
Moreover, using the results of Section 1.4 in \cite{Co}, we deduce
that $(L_S)_x$ is a maximally isotropic subspace of
$(T_xS\times\R)\oplus (T^\ast _xS\times\R)$ under the symmetric
pairing $<\, ,\, >_+$. In particular, this implies that $dim\,
(L_S)_x=dim\, S+1$, for all $x\in S$. In addition, we may prove
the following proposition.
\begin{proposition}\label{est-induc}
Let $L$ be a ${\cal E}^1(M)$-Dirac structure and $S$ be a
submanifold of $M$. If the dimension of $L_x\cap
((T_xS\times\R)\oplus (T^\ast _xM\times\R))$ keeps constant for
all $x\in S$ (or, equivalently, the dimension of $L_x\cap (\{
0\}\oplus ((T_xS)^\circ \times\{ 0\}))$ keeps constant for all
$x\in S$) then $L_S= \bigcup _{x\in S}(L_S)_x$ is a vector
sub-bundle of ${\cal E}^1(S)$ and, furthermore, $L_S$ is a ${\cal
E}^1(S)$-Dirac structure.
\end{proposition}
\prueba It is clear that $L_S$ is a maximally isotropic vector
sub-bundle of ${\cal E}^1(S)$ under the symmetric pairing $<\, ,
\, >_+$.

Now, we consider the vector bundle $\hat{L}_S$ over $S$ such that
the fiber $(\hat{L}_S)_x$ of $\hat{L}_S$ over $x\in S$ is given by
$$(\hat{L}_S)_x=L_x\cap ((T_xS\times \R)\oplus (T^\ast _xM\times
\R)).$$ Denote by $i_S:\hat{L}_S\to L$ the inclusion map, by $\pi
_S:\hat{L}_S\to L_S$ the canonical projection and by $T_L$
(respectively, $T_{L_S}$) the section of $\wedge ^3 L^\ast$
(respectively, $\wedge ^3L_S^\ast$) associated with the isotropic
vector sub-bundle $L$ (respectively, $L_S$). The map $i_S$
(respectively, $\pi _S$) is a monomorphism (respectively,
epimorphism) of vector bundles. Furthermore, if
$e_i=(X_i,f_i)+(\alpha_i,g_i)\in \Gamma (\hat{L}_S)$, with $i\in
\{ 1,2,3\}$, then, from (\ref{cuentita}) and Proposition
\ref{caracterizacion}, we get that $$
\begin{array}{clll}
(\pi _S ^\ast
T_{L_S})(e_1,e_2,e_3)&=&\displaystyle\frac{1}{2}\sum_{Cycl.(e_1,e_2,e_3)}
&\kern-60pt\Big ( i_{[X_1,X_2]}\alpha _3+g_3(X_1(f_2)-X_2(f_1)) \\
& & &\kern-60pt+X_3(i_{X_2}\alpha _1+f_2g_1)+f_3(i_{X_2}\alpha
_1+f_2g_1)\Big )\\ &=& (i_S^\ast T_L)(e_1,e_2,e_3)=0.
\end{array}
$$ Therefore, $\pi _S^\ast T_{L_S}=0$ and, since $\pi _S$ is an
epimorphism of vector bundles, we conclude that $T_{L_S}=0$. This
implies that $L_S$ is a ${\cal E}^1(S)$-Dirac structure.\QED
\subsection{The induced structure on the leaves of the characteristic
foliation of a ${\cal E}^1(M)$-Dirac
structure}\label{subprincipal} Let $L$ be a ${\cal E}^1(M)$-Dirac
structure. Denote by $(L,\makebox{{\bf [}}\, ,\, \makebox{{\bf
]}}_L,\rho _L)$ the associated Lie algebroid and by $\phi _L\in
\Gamma (L^\ast )$ the 1-cocycle defined by (\ref{1-cociclo}).

We consider the bundle map $(\rho _L,\phi _L):L\to TM\times \R$
given by
\begin{equation}\label{rho-barra}
(\rho _L,\phi _L)(e_x)=(\rho _L(e_x),\phi _L(x)(e_x)),
\end{equation}
for $e_x\in L_x$ and $x\in M$. Then, we may define the 2-form
$\Psi _L(x)$ on the vector space $(\rho _L,\phi _L)(L_x)$ by
\begin{equation}\label{2-forma-car}
\Psi _L(x)((\rho _L,\phi _L)((e_1)_x),(\rho _L,\phi
_L)((e_2)_x))=<(e_1)_x,(e_2)_x>_-,
\end{equation}
for $(e_1)_x,(e_2)_x\in L_x$, $<\, ,\, >_-$ being the natural
skew-symmetric pairing on $(T_xM\times \R)\oplus (T^\ast _xM\times
\R)$. Since $L$ is a isotropic vector sub-bundle of ${\cal
E}^1(M)$ under the symmetric pairing $<\, ,\, >_+$, we deduce that
the 2-form $\Psi _L(x)$ is well defined. Note that if $e_1,e_2\in
\Gamma (L)$ then one may consider the function $\Psi _L((\rho
_L,\phi _L) (e_1),(\rho _L,\phi _L)(e_2))\in C^\infty (M,\R)$
given by $$\Psi _L((\rho _L,\phi _L) (e_1),(\rho _L,\phi
_L)(e_2))(x)=\Psi _L(x)((\rho _L,\phi _L)((e_1)_x),(\rho _L,\phi
_L)((e_2)_x)),\mbox{ for all }x\in M.$$ In fact, if
$e_i=(X_i,f_i)+(\alpha _i,g_i)$, with $i\in \{ 1,2\}$, we have
that
\begin{equation}\label{for-dif}
\Psi _L((X_1,f_1),(X_2,f_2))=i_{X_2}\alpha _1+f_2g_1.
\end{equation}
Now, let ${\cal F}_L$ be the characteristic foliation of the
${\cal E}^1(M)$-Dirac structure $L$ and $F$ be a leaf of ${\cal
F}_L$. If $x$ is a point of $F$, we will denote by $(L_F)_x$ the
vector subspace of $(T_xF\times \R)\oplus (T^\ast _xF\times \R)$
given by (see Section \ref{induc-base}) $$(L_F)_x=\frac{L_x\cap
((T_xF\times\R)\oplus (T^\ast _xM\times\R))}{L_x\cap (\{ 0\}\oplus
((T_xF)^\circ \times\{ 0\}))}.$$ Then, we will prove that
$L_F=\bigcup _{x\in F}(L_F)_x$ defines a ${\cal E}^1(F)$-Dirac
structure and we will describe the nature of $L_F$.
\begin{theorem}\label{hojas}
Let $L$ be a ${\cal E}^1(M)$-Dirac structure and $F$ be the leaf
of the characteristic foliation ${\cal F}_L$ passing through
$x_0\in M$. Then, $L_F=\bigcup _{x\in F}(L_F)_x$ defines a ${\cal
E}^1(F)$-Dirac structure and we have two possibilities:
\begin{itemize}
\item[{\it i)}]  If $ker\,(\rho _L{ }_{|L_{x_0}})\not\subseteq <\phi
_L(x_0)>^\circ$, the ${\cal E}^1(F)$-Dirac structure $L_F$ comes
from a precontact structure $\eta _F$ on $F$, that is,
$L_F=L_{\eta _F}$. In this case, $F$ is said to be a precontact
leaf.
\item[{\it ii)}] If $ker\,(\rho _L{}_{|L_{x_0}})\subseteq <\phi
_L(x_0)>^\circ$, the ${\cal E}^1(F)$-Dirac structure $L_F$ comes
from a locally conformal presymplectic structure $(\Omega
_F,\omega _F)$ on $F$, that is, $L_F=L_{(\Omega _F,\omega _F)}$.
In this case, $F$ is said to be a locally conformal presymplectic
leaf.
\end{itemize}
\end{theorem}
\prueba From the definition of ${\cal F}_L$, it follows that
$$(\hat{L}_F)_x=L_x\cap ((T_xF\times \R)\oplus (T^\ast _xM\times
\R))=L_x,\mbox{ for all }x\in F.$$ Thus, since $L$ is a maximally
isotropic vector sub-bundle of ${\cal E}^1(M)$ under the symmetric
pairing $<\, ,\, >_+$, we obtain that $dim\, (\hat{L}_F)_x=dim\,
M+1$, for all $x\in F$. Therefore, using Proposition
\ref{est-induc}, we deduce that $L_F$ defines a ${\cal
E}^1(F)$-Dirac structure.

Next, we will distinguish the two cases:

{\it i)} Assume that $ker\,(\rho _L{ }_{|L_{x_0}})\not\subseteq
<\phi _L(x_0)>^\circ$. Then, from Theorem \ref{posibilidades}, we
have that $ker\,(\rho _L{ }_{|L_x})\not\subseteq <\phi
_L(x)>^\circ$, for all $x\in F$. This implies that the map $(\rho
_L,\phi _L)_{|L_x}:L_x\to T_xF\times \R$ is a linear epimorphism,
for all $x\in F$ (see (\ref{rho-barra})). Consequently, the
restriction of $\Psi _L$ to $F$ defines a section of the vector
bundle $\wedge ^2(T^\ast F \times\R )\to F$, i.e., a pair $(\Phi
_F,\eta _F)\in \Omega ^2(F)\times\Omega ^1(F)$. The relation
between $\Psi _L$ and $(\Phi _F,\eta _F)$ is given by
\begin{equation}\label{identifica}
\Psi _L(x)((u_1,\lambda _1),(u_2,\lambda _2))=\Phi
_F(x)(u_1,u_2)+\lambda _1\,\eta _F(x)(u_2)-\lambda _2\,\eta
_F(x)(u_1),
\end{equation}
for all $x\in F$ and $(u_1,\lambda _1),(u_2,\lambda _2)\in
T_xF\times\R$.

Now, suppose that $(u,\lambda )+(\alpha ,\mu )\in
(\hat{L}_F)_x=L_x$, with $x\in F$. From (\ref{emparejamientos}),
(\ref{2-forma-car}) and (\ref{identifica}), it follows that
$$(i_u\Phi _F(x)+\lambda\,\eta _F(x))(v)+\nu (-\eta
_F(x)(u))=\alpha (v)+\nu\,\mu,$$ for all $(v,\nu )\in T_xF\times
\R$, that is, $\alpha _{|T_xF}=i_u\Phi _F(x)+\lambda\,\eta_F(x)$
and $\mu =-\eta _F(x)(u)$. In other words, if we consider
$(L_F)_x$ to be a subspace of $(T_xF \times \R)\oplus (T^\ast
_xF\times \R)$, we have that $$(L_F)_x\subseteq\{ (u,\lambda
)+(i_u\Phi _F(x) +\lambda \,\eta _F(x),-\eta _F(x)(u))\, /
\,(u,\lambda )\in T_xF\times\R) \}.$$ But, since $dim\,
(L_F)_x=dim\, F+1$, we deduce that
 $$(L_F)_x=\{ (u,\lambda
)+(i_u\Phi _F(x) +\lambda \,\eta _F(x),-\eta _F(x)(u))\, /
\,(u,\lambda )\in T_xF\times\R) \}.$$ Thus,
\begin{equation}\label{precon}
\Gamma (L_F)=\{ (X,f )+(i_X\Phi _F +f\,\eta _F,-i_X\eta _F )\, /\,
(X,f )\in \frak X(F)\times C^\infty (F,\R) \}.
\end{equation}
Finally, using (\ref{precon}) and the fact that $L_F$ is a ${\cal
E}^1(F)$-Dirac structure, we conclude that (see Section
\ref{ejemplos}, Example 3),
\begin{equation}\label{igualdad}
\Phi _F =d\eta _F
\end{equation}
and $L_F=L_{\eta _F}$.

{\it ii)} Assume that $ker\,(\rho _L{ }_{|L_{x_0}})\subseteq <\phi
_L(x_0)>^\circ$. Then, from Theorem \ref{posibilidades}, we obtain
that $ker\,(\rho _L{ }_{|L_x})\subseteq <\phi _L(x)>^\circ$, for
all $x\in F$. This implies that the map $$\varphi
_x=pr_1{}_{|(\rho _L,\phi _L)(L_x)}:(\rho _L,\phi _L)(L_x)\to \rho
_L(L_x) =T_x F$$ is a linear isomorphism, for all $x\in F$, where
$pr_1:T_xF\times \R\to T_xF$ is the projection onto the first
factor.

Therefore, $\Psi _L$ induces a 2-form $\Omega _F$ on $F$ which is
characterized by the condition
\begin{equation}\label{2-forma}
\Omega _F(x)(\rho _L((e_1)_x),\rho _L((e_2)_x))=\Psi _L(x)((\rho
_L,\phi _L)((e_1)_x),(\rho _L,\phi _L)((e_2)_x)),
\end{equation}
for $x\in F$ and $(e_1)_x,(e_2)_x\in L_x$.

Moreover, since $S^{\phi _L}_F=F$, Theorem \ref{posibilidades}
allows us to introduce the closed 1-form $\omega _F$ on $F$
characterized by (\ref{1-form-induc}).

Now, suppose that $(u,\lambda )+(\alpha ,\mu )\in
(\hat{L}_F)_x=L_x$, with $x\in F$. From (\ref{emparejamientos}),
(\ref{1-cociclo}), (\ref{1-form-induc}), (\ref{2-forma-car}) and
(\ref{2-forma}), it follows that $\lambda =-\omega _F(x)(u)$ and
that $\alpha _{|T_xF}=i_u\Omega _F(x)+\mu\, \omega _F(x)$. In
other words, if we consider $(L_F)_x$ to be a subspace of
$(T_xF\times \R)\oplus (T^\ast _xF\times \R)$ then, since $dim\,
(L_F)_x=dim\,F+1$, we deduce that $$(L_F)_x= \{ (u,-\omega
_F(x)(u))+(i_u\Omega _F(x) +\mu \,\omega _F(x),\mu )\,/\, (u,\mu
)\in T_xF\times\R) \}.$$ Thus,
\begin{equation}\label{loc-conf}
\Gamma (L_F)=\{ (X,-\omega _F(X))+(i_X\Omega _F +f\,\omega
_F,f)\,/\,(X,f )\in \frak X(F)\times C^\infty (F,\R) \}.
\end{equation}
Finally, using that $\omega _F$ is closed, (\ref{loc-conf}) and
the fact that $L_F$ is a ${\cal E}^1(F)$-Dirac structure, we
conclude that the pair $(\Omega _F ,\omega _F)$ is a l.c.p.
structure on $F$ (see Section \ref{ejemplos}, Example 2) and that
$L_F=L_{(\Omega _F,\omega _F)}$.\QED
\begin{examples}\label{ejemp-fol}
{\rm {\bf 1.- Dirac structures}

Let $\tilde{L}\subseteq TM\oplus T^\ast M$ be a Dirac structure on
$M$ and $L$ be the ${\cal E}^1(M)$-Dirac structure associated with
$\tilde{L}$ (see (\ref{ej1})). We know that the characteristic
foliations ${\cal F}_{\tilde{L}}$ and ${\cal F}_L$ associated with
$\tilde{L}$ and $L$, respectively, coincide (see Section
\ref{ejemplos}, Example 1). Thus, if $\tilde{F}$ is a leaf of
${\cal F}_{\tilde{L}}$ then, using Theorem \ref{hojas} and the
fact that the 1-cocycle $\phi _L$ identically vanishes, it follows
that $\tilde{F}$ carries an induced l.c.p. structure $(\Omega
_{\tilde{F}},\omega _{\tilde{F}})$. Moreover, from the definition
of $\omega _{\tilde{F}}$ (see (\ref{1-form-induc})), we obtain
that $\omega _{\tilde{F}}=0$, that is, $\Omega _{\tilde{F}}$ is a
presymplectic form on $\tilde{F}$. Therefore, we deduce a
well-known result (see \cite{Co}): the leaves of the
characteristic foliation ${\cal F}_{\tilde{L}}$ are presymplectic
manifolds.

{\bf 2.- Locally conformal presymplectic structures}

Let $(\Omega ,\omega )$ be a l.c.p. structure on a manifold $M$
and $L_{(\Omega ,\omega )}$ be the corresponding ${\cal
E}^1(M)$-Dirac structure (see (\ref{ej3})). It is clear that
${\cal F}_{L_{(\Omega ,\omega)}}(x)=T_xM$, for all $x\in M$, and
thus there is only one leaf of the foliation ${\cal F}_{L_{(\Omega
,\omega)}}$, namely, $M$. Besides, since $ker\,(\rho _{L_{(\Omega
,\omega )}}{} _{|(L_{(\Omega ,\omega)})_x})\subseteq <\phi
_{L_{(\Omega ,\omega )}} (x)>^\circ$, for all $x\in M$ (see
Section \ref{ejemplos}, Example 2), $M$ carries an induced l.c.p.
structure which is just $(\Omega ,\omega )$.

{\bf 3.- Precontact structures}

Let $\eta$ be a precontact structure on a manifold $M$ and denote
by $L_\eta$ the correspon\-ding ${\cal E}^1(M)$-Dirac structure
(see (\ref{ej4})). As in the case of a l.c.p. structure, there is
only one leaf of the characteristic foliation ${\cal F}_{L_\eta}$:
the manifold $M$. In addition, since $ker\,(\rho _{L_\eta}{ }
_{|(L_\eta)_x})\not\subseteq <\phi _{L_\eta}(x)>^\circ$, for all
$x\in M$ (see Section \ref{ejemplos}, Example 3), $M$ carries an
induced precontact structure. Such a structure is defined by the
1-form $\eta$.

{\bf 4.- Jacobi structures}

Suppose that $(\Lambda ,E)$ is a Jacobi structure on a manifold
$M$ and let $L_{(\Lambda ,E)}$ be the co\-rresponding ${\cal
E}^1(M)$-Dirac structure. We know that the characteristic
foliation ${\cal F}_{L_{(\Lambda ,E)}}$ of $L_{(\Lambda ,E)}$ is
just the characteristic foliation associated with the Jacobi
structure $(\Lambda ,E)$ (see Section \ref{ejemplos}, Example 4).
Moreover, using that the Lie algebroid $(L_{(\Lambda
,E)},\makebox{{\bf [}}\, ,\, \makebox{{\bf ]}}_{L_{(\Lambda
,E)}},\rho _{L_{(\Lambda ,E)}})$ can be identified with the Lie
algebroid $(T^\ast M \times \R ,$ $\lcf \, ,\, \rcf _{(\Lambda
,E)},\widetilde{\#}_{(\Lambda ,E)})$ and that, under this
identification, the 1-cocycle $\phi _{L_{(\Lambda ,E)}}$ is the
pair $(-E,0)$, we obtain that if $x_0$ is a point of $M$ then
\begin{equation}\label{condi-Jacobi}
ker\,(\rho _{L_{(\Lambda ,E)}}{ }_{|(L_{(\Lambda ,E)})
_{x_0}})\subseteq <\phi _{L_{(\Lambda ,E)}}(x_0)>^\circ \iff
E(x_0)\in \# _\Lambda (T^\ast_{x_0} M).
\end{equation}
Thus, if $F$ is the leaf of ${\cal F}_{L_{(\Lambda ,E)}}$ passing
through $x_0\in M$ and $E(x_0)\in \# _\Lambda (T^\ast_{x_0} M)$,
from (\ref{condi-Jacobi}) and Theorem \ref{hojas}, it follows that
$F$ carries an induced l.c.p. structure $(\Omega _F,\omega _F)$.
In fact, using (\ref{emparejamientos}), (\ref{1-form-induc}),
(\ref{2-forma-car}) and (\ref{2-forma}), we have that $$\Omega
_F(y)(\#_\Lambda (\alpha _1),\# _\Lambda (\alpha _2))=\alpha _1(\#
_\Lambda (\alpha _2)),\quad \omega _F(y)(\# _\Lambda (\alpha
_1))=\alpha _1(E(y)),$$ for all $y\in F$ and $\alpha _1,\alpha_2
\in T_y^\ast M$. Therefore, the pair $(-\Omega _F,\omega _F)$ is
the locally conformal symplectic structure on $F$ induced by the
Jacobi structure $(\Lambda ,E)$.

On the other hand, if $F$ is the leaf of ${\cal F}_{L{(\Lambda
,E)}}$ passing through $x_0\in M$ and $E(x_0)\not\in \# _\Lambda
(T^\ast_{x_0} M)$ then, from (\ref{condi-Jacobi}) and Theorem
\ref{hojas}, we obtain that the ${\cal E}^1(F)$-Dirac structure
comes from a precontact structure $\eta _F$ on $F$. In addition,
$E(y)\not\in \# _\Lambda (T^\ast_y M)$ and $T_yF=\#_ \Lambda
(T^\ast _yM)\oplus $ $<E(y)>$, for all $y\in F$. Moreover, using
(\ref{emparejamientos}), (\ref{2-forma-car}) and
(\ref{identifica}), we get that $$\eta _F(y)(\#_\Lambda (\alpha
)+\lambda \, E(y))=-\lambda ,$$ for all $y\in F$, $\alpha \in
T^\ast _yM$ and $\lambda \in\R$. Consequently, $-\eta _F$ is the
contact structure on $F$ induced by the Jacobi structure $(\Lambda
,E)$.

In conclusion, we deduce a well-known result (see \cite{GL,Ki}):
the leaves of the characteristic foliation of a Jacobi manifold
are contact or locally conformal symplectic manifolds.

{\bf 5.- Homogeneous Poisson structures}

Let $(M,\Pi, Z)$ be a homogeneous Poisson manifold and $L_{(\Pi
,Z)}$ be the corresponding ${\cal E}^1(M)$-Dirac structure (see
(\ref{ej5})). Using that the Lie algebroid $(L_{(\Pi
,Z)},\makebox{{\bf [}}\, ,\, \makebox{{\bf ]}}_{L_{(\Pi ,Z)}},\rho
_{L_{(\Pi ,Z)}})$ can be identified with the Lie algebroid
$(T^\ast M \times \R ,\lcf \, ,\, \rcf _{(\Pi ,Z)},$
$\widetilde{\#}_{(\Pi ,Z)})$ and that, under this identification,
the 1-cocycle $\phi_{L_{(\Pi ,Z)}}$ is the pair $(0,1)$ (see
Section \ref{ejemplos}, Example 5), we obtain that if $x_0$ is a
point of $M$ then
\begin{equation}\label{identi}
ker\,(\rho _{L_{(\Pi ,Z)}}{ }_{|(L_{(\Pi ,Z)})_{x_0}})\subseteq
<\phi _{L_{(\Pi ,Z)}}(x_0)>^\circ \iff Z(x_0)\not\in \# _{\Pi
}(T^\ast_{x_0} M).
\end{equation}
Thus, if $x_0$ is a point of $M$ and $F$ is the leaf of the
characteristic foliation ${\cal F}_{L_{(\Pi ,Z)}}$ passing through
$x_0$, we will distinguish two cases:

{\it a)} $Z(x_0)\in \#_\Pi (T^\ast _{x_0}M)$. In such a case, from
(\ref{symp+Z}), (\ref{identi}) and Theorem \ref{posibilidades}, it
follows that $T_yF={\cal F}_{L_{(\Pi ,Z)}}(y)={\cal F}_\Pi(y)$,
for all $y\in F$, where ${\cal F}_\Pi$ is the symplectic foliation
of the Poisson manifold $(M,\Pi)$. Therefore, $F$ is the leaf of
${\cal F}_\Pi$ passing through $x_0$. In addition, using Theorem
\ref{hojas}, we deduce that the induced ${\cal E}^1(F)$-Dirac
structure comes from a precontact structure $\eta _F$ on $F$.
Moreover, from (\ref{emparejamientos}), (\ref{2-forma-car}),
(\ref{identifica}) and (\ref{igualdad}), we have that $$\eta
_F(y)(\#_\Pi (\alpha _1))=-\alpha _1(Z(y)),\quad d\eta
_F(y)(\#_\Pi (\alpha _1),\#_\Pi (\alpha _2))=\alpha _1(\#_\Pi
(\alpha _2)),$$ for all $y\in F$ and $\alpha _1,\alpha _2\in
T^\ast_yM$. This implies that $d\eta _F$ is, up to sign, the
symplectic 2-form of $F$.

{\it b)} $Z(x_0)\not\in \#_\Pi (T^\ast _{x_0}M)$. In such a case,
from (\ref{symp+Z}) and (\ref{identi}), we get that $T_yF={\cal
F}_{L_{(\Pi ,Z)}}(y)={\cal F}_\Pi(y)\oplus <Z(y)>$, for all $y\in
F$. Consequently, the dimension of $F$ is odd and the leaf $F_\Pi$
of the foliation ${\cal F}_\Pi$ passing through $x_0$ is a
submanifold of $F$ of codimension one. Furthermore, the induced
${\cal E}^1(F)$-Dirac structure comes from a l.c.p. structure
$(\Omega _F,\omega _F)$ on $F$ and, using (\ref{emparejamientos}),
(\ref{1-form-induc}), (\ref{2-forma-car}) and (\ref{2-forma}), it
follows that $$\begin{array}{c} \Omega _F(y)(\#_\Pi (\alpha
_1)+\lambda _1\,Z(y),\#_\Pi (\alpha _2)+\lambda _2\,Z(y))=\alpha
_1(\#_\Pi (\alpha _2)),\\ \omega _F(y)(\#_\Pi (\alpha _1)+\lambda
_1\,Z(y))=\lambda _1,\end{array}$$for all $y\in F$, $\alpha
_1,\alpha _2\in T^\ast _yM$ and $\lambda _1,\lambda _2\in \R$.
Note that if $i:F_\Pi\to F$ is the canonical inclusion, we deduce
that $i^\ast \omega _F=0$ and that $-i^\ast \Omega _F$ is the
symplectic 2-form on $F_\Pi$. Thus, if the dimension of $F$ is
$2n+1$, we obtain that $\omega _F\wedge \Omega _F^n=\omega
_F\wedge \Omega _F\wedge\stackrel{(n}{\ldots}\wedge \Omega _F$ is
a volume form on $F$.
 }
\end{examples}
\section{Dirac$\kern-1pt$ structure$\kern-1pt$ associated$\kern-1pt$
with$\kern-1pt$ a$\kern-1pt$ ${\cal E}^1(M)\kern-1pt
$-Dirac$\kern-1pt$ structure and characteristic foliations}
\setcounter{equation}{0} Let $M$ be a differentiable manifold and
$L$ be a vector sub-bundle of ${\cal E}^1(M)$.

We consider the vector sub-bundle $\tilde{L}$ of $T(M\times
\R)\oplus T^\ast (M\times \R)$ such that the fiber
$\tilde{L}_{(x,t)}$ of $\tilde{L}$ over $(x,t)\in M\times \R$ is
given by
\begin{equation}\label{Poissonization}
\tilde{L}_{(x,t)}=\Big \{ \Big (u+\lambda\frac{\partial }{\partial
t}_{|t}\Big )+e^t \Big ( \alpha +\mu \, dt _{|t} \Big )\, /\,
(u,\lambda )+(\alpha ,\mu )\in L_x\Big \},
\end{equation}
where $L_x$ is the fiber of $L$ over $x$. Note that the linear map
$\psi_{(x,t)}:L_x\to\tilde{L}_{(x,t)}$ given by
\begin{equation}\label{isomorfismo}
\psi _{(x,t)}((u,\lambda )+(\alpha ,\mu ))=
(u+\lambda\frac{\partial }{\partial t}_{|t} )+e^t ( \alpha +\mu \,
dt _{|t} ),
\end{equation}
is an isomorphism of vector spaces, for all $(x,t)\in M\times \R$.
Using this fact, (\ref{corchete}) and (\ref{corchete-chico}), we
deduce the following result.
\begin{proposition}\label{Diracizacion}
$L$ is a ${\cal E}^1(M)$-Dirac structure if and only if
$\tilde{L}$ is a Dirac structure on $M\times \R$.
\end{proposition}
Now, suppose that $L$ is a ${\cal E}^1(M)$-Dirac structure and
denote by $(L,\makebox{{\bf [}}\, ,\, \makebox{{\bf ]}}_L,\rho
_L)$ the associated Lie algebroid and by $\phi _L$ the 1-cocycle
of $(L,\makebox{{\bf [}}\, ,\, \makebox{{\bf ]}}_L,\rho _L)$ given
by (\ref{1-cociclo}). Then, we may consider the Lie algebroid
structure $(\makebox{{\bf [}}\, ,\, \makebox{{\bf ]}}_L
\kern-3pt\bar{ } \,^{\phi _L} ,\bar{\rho}_L^{\phi _L})$ defined by
(\ref{corchbarra}) on the vector bundle $\bar{L}=L\times\R\to
M\times\R$.

On the other hand, let $\tilde{L}$ be the Dirac structure on
$M\times \R$ associated with $L$, $(\tilde{L},\makebox{{\bf [}}\,
,\, \makebox{{\bf ]}}\,^{\tilde{ }}\kern-3pt_{\tilde{L}}
,\tilde{\rho} _{\tilde{L}})$ be the corresponding Lie algebroid
over $M\times\R$ and ${\cal F}_{\tilde{L}}$ be the characteristic
foliation of $\tilde{L}$ (see Section \ref{ejemplos}, Example 1).

It is clear that the linear maps $\psi _{(x,t)}$, $(x,t)\in
M\times \R$, induce an isomorphism of vector bundles
$\tilde{\psi}$ between $\bar{L}$ and $\tilde{L}$. Moreover, we
have
\begin{lemma}\label{lemafinal}
The map $\tilde{\psi}$ is an isomorphism of Lie algebroids over
the identity, that is,
\begin{equation}\label{isom1}
\tilde{\rho} _{\tilde{L}}(\tilde{\psi}(\bar{e}_1))=
\bar{\rho}_L^{\,\phi _L} (\bar{e}_1), \quad
\tilde{\psi}\makebox{{\bf [}}\bar{e}_1,\bar{e}_2\makebox{{\bf
]}}_L \kern-3pt\bar{ } \,^{\phi _L} =\makebox{{\bf [}}
\tilde{\psi}(\bar{e}_1) ,\tilde{\psi }(\bar{e}_2) \makebox{{\bf
]}}_{\tilde{L}}\kern-3pt\tilde{ }\kern3pt,
\end{equation}
for $\bar{e}_1,\bar{e}_2\in \Gamma (\bar{L})$. Thus, the
characteristic foliation ${\cal F}_{\tilde{L}}$ of the Dirac
structure $\tilde{L}$ coincides with the Lie algebroid foliation
${\cal F}_{\bar{L}}$.
\end{lemma}
\prueba Using (\ref{ancla}), (\ref{corchete}), (\ref{1-cociclo}),
(\ref{corchete-chico}), (\ref{ancla-chica}), (\ref{corchbarra})
and (\ref{isomorfismo}), we deduce that (\ref{isom1}) hold. In
addition, from (\ref{isom1}), it follows that ${\cal
F}_{\tilde{L}}(x,t)={\cal F}_{\bar{L}}(x,t)$, for all $(x,t)\in
M\times \R$.\QED

Now, assume that $\tilde{F}$ is a leaf of the foliation ${\cal
F}_{\tilde{L}}={\cal F}_{\bar{L}}$. Then, we know that $\tilde{F}$
is a presymplectic manifold with presymplectic 2-form
$\Omega_{\tilde{F}}$ characterized by the condition
\begin{equation}\label{2-forma-tilde}
\Omega_{\tilde{F}}(x,t)(\tilde{\rho}_{\tilde{L}}(
(\tilde{e}_1)_{(x,t)}), \tilde{\rho}
_{\tilde{L}}((\tilde{e}_2)_{(x,t)}))=<
(\tilde{e}_1)_{(x,t)},(\tilde{e}_2)_{(x,t)}>_-,
\end{equation}
for all $(x,t)\in M\times \R$ and
$(\tilde{e}_1)_{(x,t)},(\tilde{e}_2)_{(x,t)}\in
\tilde{L}_{(x,t)}$, where $<\, ,\,>_-$ is the natural
skew-symmetric pairing on $T_{(x,t)}(M\times\R)\oplus T^\ast
_{(x,t)}(M\times \R)$ (see \cite{Co} and Examples
\ref{ejemp-fol}).

Next, we will discuss the relation between the leaves of ${\cal
F}_{\tilde{L}}$ and the leaves of the characteristic foliation
${\cal F}_L$ associated with $L$. In addition, we will describe
the relation between the induced structures on them.
\begin{theorem}\label{relacion}
Let $L$ be a ${\cal E}^1(M)$-Dirac structure and $\tilde{L}$ be
the Dirac structure on $M\times\R$ associated with $L$. Suppose
that $(x_0,t_0)\in M\times \R$ and that $F$ and $\tilde{F}$ are
the leaves of ${\cal F}_L$ and ${\cal F}_{\tilde{L}}$ passing
through $x_0$ and $(x_0,t_0)$, respectively. Then:

{\it i)} If $F$ is a precontact leaf we have that
$\tilde{F}=F\times \R$. Moreover, if $\eta _F$ is the precontact
structure on $F$,  $$\Omega_{\tilde{F}}=e^t\Big ( (\pi _1{
}_{|\tilde{F}})^\ast (d\eta _F) +dt\wedge (\pi _1{
}_{|\tilde{F}})^\ast(\eta _F)\Big ),$$ where $\pi _1{
}_{|\tilde{F}}:\tilde{F}\to F$ is the restriction to $\tilde{F}$
of the canonical projection $\pi _1:M\times \R\to M$.

{\it ii)} If $F$ is a l.c.p. leaf and $(\Omega _F,\omega _F)$ is
the l.c.p. structure on $F$ then $\pi _1(\tilde{F})=F$, $\pi _1{ }
_{|\tilde{F}}: \tilde{F}\to F$ is a covering map and $\tilde{F}$
is diffeomorphic to a Galois covering of $F$ associated with
$\omega _F$. Furthermore, if $i_{\tilde{F}}:\tilde{F}\to
M\times\R$ is the canonical inclusion and $\tilde{\sigma}\in
C^\infty (\tilde{F},\R)$ is the function given by $\tilde{\sigma}
=-(i_{\tilde{F}})^\ast (t)$, we have that $$d\tilde{\sigma}=(\pi
_1{ }_{|\tilde{F}})^\ast(\omega _F),\quad
\Omega_{\tilde{F}}=e^{-\tilde{\sigma}}(\pi _1{ }
_{|\tilde{F}})^\ast(\Omega _F).$$
\end{theorem}
\prueba {\it i)} Since $F$ is a precontact leaf, it follows that
$ker\,(\rho _L{ }_{|L_{x_0}})\not\subseteq <\phi _L(x_0)>^\circ$
(see Theorem \ref{hojas}). Thus, from Theorem \ref{rel-hojas} and
Lemma \ref{lemafinal}, we deduce that $\tilde{F}=F\times\R$.

On the other hand, if $(x,t)\in \tilde{F}$,
$(\tilde{e}_i)_{(x,t)}\in \tilde{L}_{(x,t)}$, $i\in \{1,2\}$, and
$(\tilde{e}_i)_{(x,t)}=(u_i+\lambda _i\frac{\partial}{\partial
t}_{|t})+e^t(\alpha _i+\mu _idt_{|t})$, with $(u_i,\lambda
_i)+(\alpha _i,\mu _i)\in L_x$ then, using (\ref{2-forma-car}),
(\ref{identifica}), (\ref{igualdad}) and (\ref{2-forma-tilde}), we
get $$\begin{array}{rcl}\Omega_{\tilde{F}}(x,t)
(\tilde{\rho}_{\tilde{L}}((\tilde{e}_1)_{(x,t)}),
\tilde{\rho}_{\tilde{L}}((\tilde{e}_2)_{(x,t)}))&=&\frac{1}{2}e^t(\alpha
_1(u_2)+\lambda _2\mu_1-\alpha _2(u_1)-\mu_2\lambda _1)\\
\kern-200pt&\kern-200pt=&\kern-105pt e^t\Big ( (\pi _1{
}_{|\tilde{F}})^\ast (d\eta _F) +dt\wedge (\pi _1{
}_{|\tilde{F}})^\ast(\eta _F)\Big )(x,t)
(\tilde{\rho}_{\tilde{L}}((\tilde{e}_1)_{(x,t)}),
\tilde{\rho}_{\tilde{L}}((\tilde{e}_2)_{(x,t)})).
\end{array}
$$ This implies that $\Omega_{\tilde{F}}=e^t\Big ( (\pi _1{
}_{|\tilde{F}})^\ast (d\eta _F) +dt\wedge (\pi _1{
}_{|\tilde{F}})^\ast(\eta _F)\Big )$.

{\it ii)} If $F$ is a l.c.p. leaf then $ker\,(\rho _L{
}_{|L_{x_0}})\subseteq <\phi _L(x_0)>^\circ$ (see Theorem
\ref{hojas}). Therefore, from Theorem \ref{rel-hojas} and Lemma
\ref{lemafinal}, we obtain that $\pi _1(\tilde{F})=F$, that $\pi
_1{ } _{|\tilde{F}}: \tilde{F}\to F$ is a covering map, that
$\tilde{F}$ is diffeomorphic to a Galois covering of $F$
associated with $\omega _F$ and that $d\tilde{\sigma}=(\pi
_1{}_{|\tilde{F}})^\ast (\omega _F)$.

Finally, if $(x,t)\in \tilde{F}$, $(\tilde{e}_i)_{(x,t)}\in
\tilde{L}_{(x,t)}$, $i\in \{1,2\}$, and
$(\tilde{e}_i)_{(x,t)}=(u_i+\lambda _i\frac{\partial}{\partial
t}_{|t})+e^t(\alpha _i+\mu _idt_{|t})$, with $(u_i,\lambda
_i)+(\alpha _i,\mu _i)\in L_x$ then, using (\ref{2-forma-car}),
(\ref{2-forma}), (\ref{2-forma-tilde}) and the definition of
$\tilde{\sigma}$, we deduce $$\begin{array}{ccl}
\Omega_{\tilde{F}}(x,t)
(\tilde{\rho}_{\tilde{L}}((\tilde{e}_1)_{(x,t)}),
\tilde{\rho}_{\tilde{L}}((\tilde{e}_2)_{(x,t)}))&=&\kern-3pt
\frac{1}{2}e^t(\alpha _1(u_2)+\lambda _2\mu_1-\alpha
_2(u_1)-\mu_2\lambda _1)\\&=&\kern-3pt(e^{-\tilde\sigma}(\pi _1{
}_{|\tilde{F}})^\ast (\Omega
_F))(x,t)(\tilde{\rho}_{\tilde{L}}((\tilde{e}_1)_{(x,t)}),
\tilde{\rho}_{\tilde{L}}((\tilde{e}_2)_{(x,t)})).
\end{array}$$
This implies that $\Omega_{\tilde{F}}=e^{-\tilde{\sigma}}(\pi _1{
}_{|\tilde{F}})^\ast (\Omega _F)$.\QED
\begin{examples}
{\rm {\bf 1.- Dirac structures}

Let $L$ be a ${\cal E}^1(M)$-Dirac structure which comes from a
Dirac structure on $M$ and $\tilde{L}$ be the associated Dirac
structure on $M\times \R$. If $x_0$ is a point of $M$ and $F$ is
the leaf of the characteristic foliation ${\cal F}_L$ passing
through $x_0$, then $F$ is a presymplectic manifold with
presymplectic 2-form $\Omega _F$ (see Examples \ref{ejemp-fol}).
Moreover, since ${\cal F}_L(x)={\cal F}_{\tilde{L}}(x,t)$, for all
$(x,t)\in M\times \R$, we deduce that the leaf $\tilde{F}$ of
${\cal F}_{\tilde{L}}$ passing through $(x_0,t_0)\in M\times \R$
is $\tilde{F}=F\times \{t_0\}$. In addition, from Theorem
\ref{relacion}, it follows that $\Omega_{\tilde{F}}=e^{t_0}\Omega
_F$.

{\bf 2.- Precontact structures}

Let $\eta$ be a precontact structure on a manifold $M$ and
$L_\eta$ be the associated ${\cal E}^1(M)$-Dirac structure. Then,
the characteristic foliation ${\cal F}_{L_\eta}$ has a unique
leaf, the manifold $M$ (see Examples \ref{ejemp-fol}).
Furthermore, if $\tilde{L}_\eta$ is the Dirac structure on
$M\times\R$ associated with $L_\eta$, we obtain that
$\tilde{L}_\eta$ is the graph of the presymplectic 2-form
$\tilde{\Omega}$ on $M\times\R$ given by $$\tilde{\Omega}=e^t((\pi
_1)^\ast (d\eta) +d t\wedge (\pi _1)^\ast \eta ).$$ In other
words, $$\Gamma (\tilde{L}_\eta)=\{ \tilde{X}+i_{\tilde{X}}
\tilde{\Omega}\, /\, \tilde{X}\in \frak X(M\times \R)\}\subseteq
\frak X (M\times \R)\oplus \Omega ^1(M\times \R) .$$ On the other
hand, using Theorem \ref{relacion}, we deduce a well-known result
(see \cite{Co}): the characteristic foliation ${\cal
F}_{\tilde{L}_\eta}$ of $\tilde{L}_\eta$ also has a unique leaf
$\tilde{F}$ (the manifold $M\times\R$) and the presymplectic
2-form $\Omega _{\tilde{F}}$ on $\tilde{F}$ is just
$\tilde{\Omega}$.

{\bf 3.- Jacobi structures}

Suppose that $(\Lambda ,E)$ is a Jacobi structure on a manifold
$M$. Then, it is well-known that the 2-vector $\tilde{\Lambda}$ on
$M\times\R$ given by $\tilde{\Lambda}=e^{-t}(\Lambda
+\frac{\partial}{\partial t}\wedge E)$ defines a Poisson structure
on $M\times\R$ (see \cite{Li2}). Thus, one may consider the Dirac
structure $\tilde{L}_{\tilde{\Lambda}}$ on $M\times\R$ associated
with $\tilde{\Lambda}$ (see \cite{Co}). In fact, we have that $$
\Gamma (\tilde{L}_{\tilde{\Lambda}})=\{
\#_{\tilde{\Lambda}}(\tilde{\alpha} )+\tilde{\alpha }\, /\,
\tilde{\alpha} \in \Omega ^1 (M\times \R)\}.$$ Moreover, if
$L_{(\Lambda ,E)}$ is the ${\cal E}^1(M)$-Dirac structure induced
by the Jacobi structure $(\Lambda ,E)$, it is easy to prove that
the Dirac structure $\tilde{L}_{(\Lambda ,E)}$ on $M\times\R$
associated with $L_{(\Lambda ,E)}$ is isomorphic to
$\tilde{L}_{\tilde{\Lambda}}$. Therefore, using Theorem
\ref{relacion} (see also Examples \ref{ejemp-fol}), we directly
deduce the results of Guedira-Lichnerowicz (see Section III.16 in
\cite{GL}) about the relation between the leaves of the
characteristic foliation of the Jacobi manifold $(M,\Lambda ,E)$
and the leaves of the symplectic foliation of the Poisson manifold
$(M\times\R ,\widetilde{\Lambda})$. }
\end{examples}

\vspace{.2cm} {\small {\bf Acknowledgments.} Research partially
supported by DGICYT grant BFM2000-0808 (Spain). D. Iglesias wishes
to thank Spanish Ministerio de Educaci\'on y Cultura for a FPU
grant.}


\vspace{-.3cm}
\begin{thebibliography}{9}
{\small
\bibitem{BV} K.H. Bhaskara, K.Viswanath: {\em Poisson algebras
and Poisson manifolds}, Research Notes in Mathema\-tics, 174,
Pitman, London, 1988.

\vspace{-9pt}
\bibitem{Co} T.J. Courant: Dirac manifolds, {\em Trans. A.M.S.} {\bf
319} (1990), 631-661.

\vspace{-9pt}
\bibitem{CW} T.J. Courant, A. Weinstein: Beyond Poisson
structures, In s\'eminaire Sud-Rhodanien de geom\'etrie, Travaux
en Cours, {\bf 27}, Hermann, Paris, 1988, 39-49.

\vspace{-9pt}
\bibitem{DLM} P. Dazord, A. Lichnerowicz and Ch.M. Marle: Structure
locale des vari\'et\'es de Jacobi, {\em J. Math. Pures Appl.},
{\bf 70} (1991), 101-152.

\vspace{-9pt}
\bibitem{Di} P.A.M. Dirac: {\em Lectures on Quamtum Mechanics},
Yeshiva Univ., New York, 1964.

\vspace{-9pt}
\bibitem{Do} I. Dorfman: {\em Dirac structures and integrability
of nonlinear evolution equations}, Nonlinear Science: Theory and
Applications. John Wiley \& Sons, Ltd., Chichester, 1993.

\vspace{-9pt}
\bibitem{Go} C. Godbillon: {\em \'{E}l\'{e}ments de topologie alg\'{e}brique},
Hermann, Paris, 1971.

\vspace{-9pt}
\bibitem{GHV} W. Greub, S. Halperin, R. Vanstone: {\em
Connections, curvature and cohomology. Vol. I}, Pure and Applied
Mathematics 47, Academic Press, New York-London, 1972.

\vspace{-9pt}
\bibitem{GL} F. Gu\'edira, A. Lichnerowicz: G\'eom\'etries des
alg\'ebres de Lie locales de Kirillov, {\em J. Math. Pures Appl.},
{\bf 63} (1984), 407-484.

\vspace{-9pt}
\bibitem{HM} P.J. Higgins, K.C.H. Mackenzie: Algebraic constructions
in the category of Lie algebroids, {\em J. Algebra}, {\bf 129}
(1990), 194-230.

\vspace{-9pt}
\bibitem{IM} D. Iglesias, J.C. Marrero: Some linear Jacobi
structures on vector bundles, {\em C.R. Acad. Sci. Paris}, {\bf
331} S\'er. I (2000), 125-130. {\em arXiv: math.DG/0007138}.

\vspace{-9pt}
\bibitem{IM2} D. Iglesias, J.C. Marrero: Generalized Lie bialgebroids
and Jacobi structures, {\em To appear in J. Geom. Phys.} (2001)
{\em arXiv: math.DG/0008105}.

\vspace{-9pt}
\bibitem{IM3} D. Iglesias, J.C. Marrero: Generalized Lie
bialgebras and Jacobi structures on Lie groups, {\em Preprint
(2001), arXiv: math.DG/0102171}.

\vspace{-9pt}
\bibitem{KS} Y. Kerbrat, Z.Souici-Benhammadi: Vari\'{e}t\'{e}s de Jacobi et
groupo\"{\i}des de contact, {\em C.R. Acad. Sci. Paris}, {\bf 317} S\'{e}r.
I (1993), 81-86.

\vspace{-9pt}
\bibitem{Ki} A. Kirillov: Local Lie algebras, {\em Russian Math.
Surveys}, {\bf 31} (1976), 55-75.

\vspace{-9pt}
\bibitem{LM} P. Libermann, Ch. M. Marle: {\em Symplectic Geometry and
Analytical Mechanics}, Kluwer, Dordrecht, 1987.

\vspace{-9pt}
\bibitem{Li1} A. Lichnerowicz: Les vari\'et\'es de Poisson et leurs
alg\'ebres de Lie associ\'ees, {\em J. Differential Geometry},
{\bf 12} (1977), 253-300.

\vspace{-9pt}
\bibitem{Li2} A. Lichnerowicz: Les vari\'et\'es de Jacobi  et leurs
alg\'ebres de Lie associ\'ees, {\em J. Math. Pures Appl.}, {\bf
57} (1978), 453-488.

\vspace{-9pt}
\bibitem{LWX} Z-J. Liu, A. Weinstein, P. Xu: Manin triples for Lie
bialgebroids, {\em J. Differential Geom.}, {\bf 45} (1997),
547-574.

\vspace{-9pt}
\bibitem{LWX2} Z-J. Liu, A. Weinstein, P. Xu: Dirac structures and
Poisson homogeneous spaces, {\em Comm. Math. Phys.}, {\bf 192}
(1998), 121-144.

\vspace{-9pt}
\bibitem{Mk} K. Mackenzie: {\em Lie groupoids and Lie algebroids in
differential geometry}, Cambridge University Press, 1987.

\vspace{-9pt}
\bibitem{MX} K. Mackenzie, P. Xu: Lie bialgebroids and Poisson
groupoids, {\em Duke Math. J.}, {\bf 73} (1994), 415-452.

\vspace{-9pt}
\bibitem{Mo} T. Mokri: Matched pairs of Lie
algebroids, {\em Glasgow Math. J.}, {\bf 39} (1997), 167-181.

\vspace{-9pt}
\bibitem{NVQ} Ng\^o-van-Qu\^e: Sur l'espace de prolongement
diff\'erentiable, {\em J. Differential Geometry}, {\bf 2} (1968),
33-40.

\vspace{-9pt}
\bibitem{RW} D. Roytenberg, A. Weinstein: Courant algebroids and
strongly homotopy Lie algebras, {\em Lett. Math. Phys.}, {\bf 46}
(1998), 81-93.

\vspace{-9pt}
\bibitem{Su} H. Sussman: Orbits of families of vector fields and
integrability of distributions, {\em Trans. A.M.S.} {\bf 180}
(1973), 171-188.

\vspace{-9pt}
\bibitem{V0} I. Vaisman: Locally conformal symplectic manifolds, {\em
Internat. J. Math. \& Math. Sci.}, {\bf 8} (1985), 521-536.

\vspace{-9pt}
\bibitem{V} I. Vaisman: {\em Lectures on the Geometry of Poisson
Manifolds}, Progress in Math. 118, Birkh\"auser, Basel, 1994.

\vspace{-9pt}
\bibitem{V2} I. Vaisman: The BV-algebra of a Jacobi manifold, {\em
Ann. Polon. Math.}, {\bf 73} (2000), 275-290. {\em arXiv: math.DG/
9904112}.

\vspace{-9pt}
\bibitem{Wa} A. Wade: Conformal Dirac Structures, {\em Lett. Math.
Phys.}, {\bf 53} (2000), 331-348. {\em arXiv: math.SG/ 0101181}.

\vspace{-9pt}
\bibitem{We} A. Weinstein: The local structure of Poisson manifolds,
{\em J. Differential Geometry}, {\bf 18} (1983), 523-557, Errata
et Addenda {\bf 22} (1985), 255. }
\end{thebibliography}
\end{document}